\documentclass[graybox]{svmult}

\usepackage[show]{ed}
\usepackage{multirow}
\usepackage{amsmath}
\usepackage{amsfonts}
\usepackage{amssymb}%
\usepackage{latexsym}
\usepackage[dvips]{graphicx}
\usepackage{psfrag}
\usepackage{color} 

\usepackage{mathptmx}       % selects Times Roman as basic font
\usepackage{helvet}         % selects Helvetica as sans-serif font
\usepackage{courier}        % selects Courier as typewriter font
\usepackage{type1cm}        % activate if the above 3 fonts are
\usepackage{makeidx}         % allows index generation
\usepackage{graphicx}        % standard LaTeX graphics tool
\usepackage{multicol}        % used for the two-column index
\usepackage[bottom]{footmisc}% places footnotes at page bottom

\makeindex             % used for the subject index
\newcommand{\Hsub}{H_{\mathrm{sub}}}
\newcommand{\ksqeps}{\frac{k^2}{ \eps}}
\newcommand{\bfv}{\mathbf{v}}
\newcommand{\bfx}{\mathbf{x}}
\newcommand{\bfd}{\mathbf{d}}
\newcommand{\bfr}{\mathbf{r}}
\newcommand{\br}{\mathbf{r}}
\newcommand{\bv}{\mathbf{v}}
\newcommand{\bw}{\mathbf{w}}
\newcommand{\C}{\mathbb{C}}
\newcommand{\bW}{\mathbf{W}}
\newcommand{\matrixR}{{ R}}
\newcommand{\cK}{\mathcal{K}}
\newcommand{\bV}{\mathbf{V}}
\newcommand{\matrixI}{{ I}}

\newcommand{\matrixP}{{ P}}

\newcommand{\matrixB}{{ B}}
\newcommand{\matrixC}{{ C}}

\newcommand{\matrixBepsRAS}{{ B}_{\eps,RAS}^{-1}}
\newcommand{\matrixA}{A}
\newcommand{\matrixAeps}{{A}_\eps}
\newcommand{\matrixAepsinv}{\matrixAeps^{-1}}

\newcommand{\cO}{\mathcal{O}}
\newcommand{\eps}{\varepsilon}
\newcommand{\epsprob}{\varepsilon_{\mathrm{prob}}}
\newcommand{\epsprec}{\varepsilon_{\mathrm{prec}}}
\newcommand{\hprob}{h_{\mathrm{prob}}}
\newcommand{\Hprec}{H_{\mathrm{prec}}}
\newcommand{\rhoprob}{\rho_{\mathrm{prob}}}
\newcommand{\rhoprec}{\rho_{\mathrm{prec}}}
\newcommand{\ri}{\mathrm{i}}

\newcommand{\cI}{{\mathcal I}}

\begin{document}

\title*{Recent Results on Domain Decomposition Preconditioning for the High-frequency Helmholtz Equation using   
Absorption }
\titlerunning{Domain Decomposition Preconditioning for High-frequency  Helmholtz} 
%for an abbreviated version of
% your contribution title if the original one is too long
\author{I.G. Graham,  E. A. Spence and E. Vainikko }
% Use \authorrunning{Short Title} for an abbreviated version of
% your contribution title if the original one is too long
\institute{Ivan G. Graham  \at Department of Mathematical Sciences, University of Bath, Bath BA2 7AY, UK \email{I.G.Graham@bath.ac.uk}
\and Euan A. Spence  \at   Department of Mathematical Sciences, University of Bath, Bath BA2 7AY, UK \email{E.A.Spence@bath.ac.uk}\and Eero Vainikko \at Institute of Computer Science, University of Tartu, Tartu 50409, Estonia \email{eero.vainikko@ut.ee} }
%
% Use the package "url.sty" to avoid
% problems with special characters
% used in your e-mail or web address
%
\maketitle

To appear in: 
{\em Modern Solvers for Helmholtz problems}\ 
 \ edited by {\sc D. Lahaye, J. Tang and C. Vuik}, 
\quad Springer Geosystems Mathematics series, 
2016.  

\vspace{0.7in} 

\abstract{In this paper we present  an overview of recent progress on the 
development and analysis of domain decomposition preconditioners for 
discretised  Helmholtz problems, where the preconditioner is constructed from 
the corresponding problem with added   absorption. Our  preconditioners incorporate local 
subproblems that can have various  boundary 
conditions,  and include the possibility of a global coarse mesh.    
{While the rigorous analysis describes preconditioners for the Helmholtz  
problem with added absorption,  this theory also informs the  development of    
efficient  multilevel solvers for the ``pure''  Helmholtz problem without absorption. 
For this case, 2D experiments for problems containing up to about  $50$ wavelengths are presented.} The {experiments} show iteration counts of order about  $\mathcal{O}(n^{0.2})$ and   
times (on a serial machine)  of order about  
 $\mathcal{O}(n^{\alpha})$, { with $\alpha \in [1.3,1.4]$}  for solving systems of dimension  $n$. This holds both in the   
pollution-free case   corresponding to meshes with 
grid size   $\mathcal{O}(k^{-3/2})$ (as the wavenumber $k$ increases), 
and also for discretisations with a fixed number of grid points per wavelength, commonly used in applications. 
% For  systems arising from a fixed number of grid-points per wavelength, 
% $n \sim k^2$ and the time for computation is close to $\mathcal{O}(n)$.  }
{Parallelisation of the algorithms is also  briefly discussed.} }
\section{Introduction}
\label{sec:Intro}

In this paper we describe recent work on the theory and implementation  of domain decomposition methods 
for iterative solution of   discretisations of the  Helmholtz equation:   
\begin{equation}\label{eq:PDE0} 
-(\Delta   + k^2)u = f \ , \quad \text{in a domain} \quad \Omega ,
\end{equation} 
  where  $k(x) = \omega/c(x)$,  with $\omega$ denoting frequency and $c$ denoting the speed of acoustic waves 
in $\Omega$. Our motivation originates from applications in seismic imaging, but the 
methods developed are applicable more generally, e.g. to earthquake modelling or medical imaging.  
While practical imaging problems often  involve the frequency domain reduction of the elastic wave 
equation or Maxwell's equations, the scalar Helmholtz equation \eqref{eq:PDE0} is still an 
extremely relevant model problem which encapsulates many of the key difficulties of more 
complex problems. 

We will focus here on solving  \eqref{eq:PDE} on a bounded domain $\Omega$,  subject to the first order absorbing   
(impedance) boundary condition:
\begin{equation}\label{eq:ImpBC} 
\frac{\partial u }{\partial n} - i k u = g \quad \text{on} \quad \Gamma = \partial \Omega \ , 
\end{equation}
although the methods presented are more general. 

The theoretical part of this paper is restricted to the case of $k$ constant. However the methods proposed can be used in the variable $k$ case, and preliminary experiments are done on this case in \S \ref{subsec:Variable}. 

Important background  for our investigation  is  the  large  body of work on ``shifted Laplace'' 
preconditioning  for this problem, starting from \cite{ErVuOo:04} and including,  
for example \cite{ErOoVu:06} and recent {work} on deflation \cite{ShLaVu:13}.  (A fuller survey is given in \cite{GaGrSp:15, GrSpVa:15} and elsewhere in this volume.) 
In those papers     (multigrid) 
approximations  of the solution operator for the 
perturbed problem  
\begin{equation}\label{eq:PDE} 
-(\Delta   + (k^2+ \ri \eps) )u = f \ , \quad \text{with} \quad \frac{\partial u }{\partial n} - i k u = g \quad \text{on} \quad \Gamma  \ , 
\end{equation} 
(suitably discretised and with carefully tuned   ``absorption''  parameter $\eps > 0 $),  
were used as preconditioners for the iterative solution of  \eqref{eq:PDE0}. 
When $k$ is variable, a slightly different shift strategy is appropriate (see \S \ref{subsec:Variable}). 

One can see immediately the benefit of introducing $\eps$ in \eqref{eq:PDE}:  
When $k$ is constant the fundamental solution  $G_{k,\eps}$ 
of the operator in \eqref{eq:PDE} (for example in 3D) satisfies,  for fixed $x\not = y$ 
% and  fixed 
% $\eps>0$, $$ G_{k,\eps}(x,y)  \sim G_{k,0}(x,y) \exp\left( -\frac{\eps}{k} \vert x - y\vert \right)\ , \quad \text{as} \quad   k \rightarrow \infty. $$     \ednote{Euan please check}
with $k\vert x-y \vert = \mathcal{O}(1)$ and $\eps\ll k^2$,  
$$ G_{k,\eps}(x,y)  =G_{k,0}(x,y) \exp\left( -\frac{\eps}{2k} \vert x - y\vert \right)\left(1 + \mathcal{O}\left( \left( \frac{\eps}{k^2}\right)^2 k\vert x-y\vert\right)\right)
\ , \quad \text{as} \quad   k \rightarrow \infty. $$  

Thus,   the effect of introducing   
$\eps$ is to exponentially damp the oscillations in the fundamental solution of problem \eqref{eq:PDE}, with the amount of damping proportional to $\eps/k$.  With slightly more analysis one can show  
that the weak form of problem 
\eqref{eq:PDE}   enjoys a coercivity property (with coercivity 
constant of order $\mathcal{O}(\eps/k^2)$  {in the energy norm \eqref{eq:energy}}   \cite[Lemma 2.4]{GrSpVa:15}). This  has the useful ramification that 
any finite element method for   \eqref{eq:PDE}  is always well-posed (independent of mesh size) and enjoys a corresponding 
(albeit $\eps-$ and $k-$ dependent)  quasioptimality property.  Therefore preconditioners constructed by applying  local  
and coarse mesh solves applied to \eqref{eq:PDE}  are always well-defined;  this is not true when $\eps = 0$.  

A natural question is then, how should one choose $\eps$? To begin to investigate this question, we first introduce some notation.
Let  $A_\eps$ denote the finite element approximation of \eqref{eq:PDE}  and write $A = A_0$. Then $A$   
 is the system matrix for problem \eqref{eq:PDE0}, \eqref{eq:ImpBC},     which  we want to solve. 

Suppose an approximate inverse  $B_\eps^{-1}$ for   $A_\eps$ is constructed. Then a sufficient condition for 
$B_\eps^{-1}$ to be a good preconditioner for $A$ is that 
$I - B_\eps^{-1}A$ should be sufficiently small. Writing
$$ I - B_\eps^{-1}A = (I - B_\eps^{-1} A_\eps) + B_\eps^{-1}A_\eps ( I - A_\eps^{-1} A ),  $$ 
% $$ \matrixI -\matrixBepsinv\matrixA \ = \ \matrixI - \matrixBepsinv \matrixA_\eps
% \ + \ \matrixBepsinv \matrixA_\eps (\matrixI - \matrixA_\eps^{-1} \matrixA) , 
% %(\bI- \matrixAeps^{-1}\matrixA) + (\tildematrixAepsinv- \matrixAepsinv)\matrixA  
% $$
% $$  \Vert I - B_\eps^{-1} A \Vert_2 \
% \text{ and }  
%  \ \  \Vert \matrixI -
% \matrixBepsinv \matrixAeps  \Vert_2  \quad \text{small} \ ,$$
we see that a sufficient condition for the smallness of the term on the left-hand side is that 
{ \begin{itemize}
\item[(i)] \ $I - A_\eps^{-1}A$ should be sufficiently small,  and 
\item[(ii)] \  $I - B_\eps^{-1}A_\eps$ should be sufficiently small.   
\end{itemize}}
At this stage, one might already guess that achieving both (i) and (ii) imposes somewhat contradictory requirements on $\eps$. Indeed, on the one hand, (i) requires $\eps$ to be sufficiently small (since the ideal preconditioner for $\matrixA$ is $\matrixA^{-1} \ = \ \matrixA_0^{-1}$). On the other hand, the larger $\eps$ is, the less oscillatory the shifted problem is, and the easier  it should  be to construct a good approximation to $\matrixAepsinv$ for (ii). 

Regarding (i): The spectral analysis in \cite{ErGa:12} of a 1-d finite-difference discretisation concluded that one needs 
 $\eps<k$  for the eigenvalues to be clustered around $1$ (which partially achieves (i)). The analysis in \cite{GaGrSp:15} showed that, in both 2- and 3-d  for a range of geometries and  
 finite element discretisations,  (i) is guaranteed if $\eps/k \leq C_1 $  for a small enough positive constant $C_1$, with numerical experiments indicating that this condition is sharp.
Somewhat different investigations are   contained in the references 
\cite{ErVuOo:04}, \cite{Er:08}, \cite{ErOoVu:06}.  These performed 
spectral analyses that essentially aim to 
achieve (i) on a continuous level, and explored the best preconditioner  
of the form \eqref{eq:PDE} for \eqref{eq:PDE0} in the 1D case with Dirichlet boundary conditions,   
based on the ansatz $k^2 + \ri \eps = k^2(a + \ri b)$, where $a,b$ are to be chosen; related more general results are in \cite{VaErVu:07}.   
(For more detail, see, e.g., the summary in \cite{GaGrSp:15} and other articles in this volume.)

Regarding (ii): several authors have considered the question of when multigrid converges (in a $k$-independent number of steps) when applied to the shifted problem $A_\eps$, with the conclusion that one needs $\eps \sim k^2$  \cite{CoGa:14}, \cite{CoVa:13}, \cite{ErGa:12}. Note that this question of convergence is not quite the same question as whether  a multigrid approximation to $A_{\eps}^{-1}$ is a good preconditioner for $A_\eps$ (property (ii)) or for $A_0$ (the original problem), but these questions are investigated numerically in \cite{CoGa:14}.
For classical Additive Schwarz 
domain decomposition preconditioners, it was shown in \cite{GrSpVa:15} that (ii) is guaranteed (under certain conditions on the coarse grid diameter) 
if $\eps\sim k^2$ (resonating with the multigrid results).
%if $\eps/k^2  \geq  C_2$ for any positive constant $C_2$.  
In fact  \cite{GrSpVa:15} also provides $\eps$-explicit estimates of the rate of GMRES convergence 
when $A_\eps$ is preconditioned by the Schwarz algorithm. Although  
these estimates degrade sharply when $\eps$ is chosen less than $k^2$,  numerical experiments in \cite{GrSpVa:15}  indicate that improved estimates may be possible in the range  $k\lesssim \eps \lesssim  k^2$.

The contradictory requirements that (i) requires $\eps/k$ to be sufficiently small, and (ii) requires $\eps\sim k^2$ (at least for classical Additive Schwarz 
domain decomposition preconditioners)  motivate the question of whether 
new choices of $B_\eps^{-1}$ 
can be devised that operate best when $\eps$ is chosen in the range  { $k \lesssim   \eps \lesssim  k^2$}. 
Such choices should necessarily use components that are more suitable for ``wave-like'' problems, rather than the essentially ``elliptic'' technology of classical multigrid or classical domain decomposition. 
In fact our numerical experiments below indicate that, for the preconditioners studied here, the best choice of $\eps$ varies,  but is generally in the range $[k, k^{1.6}]$.

Domain decomposition methods offer the attractive feature that their coarse grid and local 
problems can be adapted to allow for ``wave-like'' behaviour. There is indeed a large 
literature on this (e.g. \cite{BeDe:97,Fa:00,GaMaNa:02}, but  methods that combine  
 many subdomains and  coarse grids and include    a convergence analysis  are still missing.        
 The paper \cite{GrSpVa:15} 
provides the first such rigorous analysis in the many subdomain case, and current work is focused on extending this to the case when wave-like components are inserted, such as using  (optimised) impedance or PML conditions on 
the local solves.       
 
Another class of preconditioners for Helmholtz problems 
of great recent interest  is the ``sweeping'' preconditioner 
\cite{EnYi10} and its related variants - e.g. \cite{Chen13a}, \cite{Stolk13}, 
{\cite{ZeDe:14}. In principle these methods require 
the direct solution of Helmholtz subproblems on strips of the domain. 
A method of expediting these inner solves with an additional domain decomposition and off-line computation of local inverses is presented in \cite{ZeDe:15}.}  
Related domain decomposition methods for these inner solves, using tuned  absorption,  and  with applications to industrial problems, 
 are explored in 
\cite{ChGrSh:13}, \cite{Sh:15},    \cite{ShChGr:16}. 
 
Finally it should be acknowledged that,  while 
the  reduction of the  complicated question of the performance of $B_\eps^{-1} $ as a 
preconditioner for $A$ into two digestible subproblems  ((i) and (ii) above)  
is theoretically convenient, this approach is  also  very crude in several ways: 
Firstly the splitting of the problem into   (i) and (ii) may 
not be optimal and secondly the overarching requirement that $\Vert I - B_\eps^{-1}A\Vert $ should be  small 
is far from necessary when assessing $B_\eps^{-1} $ as a preconditioner for $A$: 
for example  good GMRES convergence is  still assured  if the field of values 
of $B_\eps^{-1} A$    is bounded away from the origin in the complex plane 
(in a suitable inner product) and that $B_\eps^{-1}A $ 
is bounded from above in the corresponding norm.    We use this in the theory below.

\section{Domain Decomposition} 
\label{sec:DD} 

To start, we denote the nodes of the  finite element mesh as $\{x_j: j \in \cI^h\}$,  for a suitable index set  
 $\cI^h$. These include  nodes on the boundary $\Gamma $ of $\Omega$.   {The continuous piecewise linear finite element hat function basis is denoted     $\{ \phi_j : j \in \cI^h\}$.}
To define preconditioners,   we  choose   a collection of $N$ non-empty  relatively open subsets 
$\Omega_\ell$  of $\overline{\Omega}$,  which    
form an overlapping cover of $\overline{\Omega}$. 
Each $\overline{\Omega}_\ell$  is assumed to  consist of a union of
elements of the finite element mesh, and the corresponding nodes on $\Omega_\ell$ are denoted $\{x_j: j \in \cI^h(\Omega_\ell)\}$. 
 % Then, for each $\ell$,  we set 
% $$\cV_\ell = \{ v_h \in \cV^h: \supp(v_h) \subset
% \overline{\Omega}_\ell\}. $$  Note that, since functions in $\cV^h$
% are continuous, functions in  $\cV_\ell$ must vanish  on the
% internal boundary $\partial \Omega_\ell \backslash \Gamma$,
% but are unconstrained on the
% external boundary $ \partial \Omega_\ell\cap \Gamma$.
% The freedoms for  $\cV_\ell$ are denoted   
% $\cN^h(\Omega_\ell) = \{x_j:j \in \cI^h(\Omega_\ell)\}$, where $\cI^h(\Omega_\ell)$ is a suitable
% index set. The  
% basis for
% $\cV^h(\Omega_\ell)$ can then be  written   $\{ \phi_j : j \in \cI^h(\Omega_\ell)
% \}$. 
% Thus a solve of the Helmholtz problem \eqref{eq:vf_intro}  in the space $\cV_\ell$ involves a
% Dirichlet 
% boundary condition at internal boundaries and natural boundary
% condition at external boundaries (if any). The introduction of the absorption
% $\eps \not = 0$ ensures such solves are always  well-defined.  Future work
% will consider  the analysis of methods with other local boundary
% conditions (such as impedance or PML).  Internal impedance
% conditions are considered in the experiments in \S \ref{sec:Numerical}.  

Now, for  any  $j \in \cI^h({\Omega_\ell})$ and $j' \in \cI^h$, we define the
restriction matrix  
$(R_\ell)_{j,j'} := \delta_{j,j'}$. The  matrix 
$$A_{\eps,\ell} \ := \ R_\ell  A_\eps R_\ell ^T$$
is then  just the minor  of 
$A_\eps$ corresponding to rows
and columns taken from $\cI^h({\Omega_\ell})$.  This matrix corresponds to a discretisation (on the fine mesh) 
of the original problem \eqref{eq:PDE} restricted to the local domain  $\Omega_\ell$, with a 
homogeneous Dirichlet condition at the 
interior boundary $\partial \Omega_\ell \backslash \Gamma$ and impedance condition 
at the 
outer boundary  $\partial \Omega_\ell \cap  \Gamma$  (when this is non-empty).    

One-level domain
decomposition methods are constructed from the inverses $A_{\eps,\ell}^{-1}$. 
More precisely, 
\begin{equation}\label{eq:minor} 
 B_{\eps, {AS}, local}^{-1} \ : = \  \sum_{\ell}R_\ell^T A_{\eps,\ell}^{-1} R_\ell,   
\end{equation}
 is the classical
one-level  Additive Schwarz approximation of  $A_\eps^{-1}$ with the subscript  ``$local$''
indicating  that the solves are on local subdomains $\Omega_\ell$. 

The overlapping subdomains are required to satisfy certain technical conditions concerning their shape and the 
size and uniformity of the overlap. Moreover,  each point in the domain is allowed to  lie 
only in a bounded  number of overlapping subdomains as the mesh is refined. 
We do not repeat these conditions  
here but refer the interested reader to \cite[\S 3]{GrSpVa:15}.   
The theorems presented in \S \ref{sec:Main} require these 
assumptions for their proof, as well as a quasi-uniformity assumption on the coarse mesh 
which is introduced next. 

% For the theory, we need  assumptions on the shape of the subdomains and
% the size of the overlap, and we require any point in $\overline{\Omega}$ to  belong to a
% bounded number of overlapping subdomains.   

Two-level methods are obtained by adding a  global coarse solve.  
% Let us introduce a family of  shape-regular, simplicial meshes on 
% $\overline{\Omega}$. 
We introduce a family of coarse simplicial  meshes with nodes $\{x_j^H,   j \in \cI^H\}$,  where each  coarse element 
is also assumed to  consist of the union of a set of 
fine grid elements. The  basis functions  are taken to be the continuous $P_1$ hat 
functions on the coarse mesh, which we denote $\{\Phi^H_p, \ p \in \cI^H\}$. 
% From these functions we define 
% the coarse space
% $
% \cV_0 \ := \ \mathrm{span} \{ \Phi_p : p \in \cI^{H}
% \}\ , 
% $
% which is  a subspace
% of $\cV^h$. 
Then,  introducing  the fine-to-coarse restriction matrix 
$(R_0)_{pj} \ := \ \Phi_p^H(x_j^h)\ ,  \ j \in \cI^h , \
p \in \cI^{H},$ 
we can define the corresponding coarse mesh matrix 
%\begin{equation}\label{eq:coarsegrid}
$A_{\eps,0} := R_0 A_\eps R_0^T\ $ . 
%\end{equation}
% This is  the stiffness matrix for problem \eqref{eq:PDE},  
% discretised in the coarse space using the basis $\{\Phi_p^H: p \in
% \cI^{H}\}$.  
Note that, due to the coercivity property for problem \eqref{eq:PDE},   
both $A_{\eps,0}$ and $A_{\eps,\ell}$ are invertible for all mesh 
sizes $h, H$ and all choices of  $\epsilon \not = 0$. 

The classical Additive 
Schwarz preconditioner  is then 
\begin{equation}\label{eq:defAS}
B_{\eps,AS}^{-1} \ : = \  R_0^T A_{\eps,0}^{-1} R_0 \ + \  B_{\eps,AS,local}^{-1} , 
\end{equation}
(i.e.~the sum of  coarse solve and local solves) with
$B_{\eps,AS,local}^{-1}$ defined in \eqref{eq:minor}.  

The theoretical results outlined in the next section concern the properties of $B_{\eps, AS}^{-1}$  
as a preconditioner for $A_\eps$ (i.e. criterion (ii)  in \S \ref{sec:Intro}).  
The hypotheses for the theory involve conditions on  $k$, $\eps$ and $H$ (the coarse mesh diameter) as well as  {$\Hsub$} (the maximum of the diameters of   the local subdomains $\Omega_\ell$).    This theory is verified by some of the   numerical experiments in \cite{GrSpVa:15} and we do not repeat those here.   
Instead, in  \S \ref{sec:Experiments} below we  focus in detail on the performance of  (variants of) 
$B_{\eps, AS}^{-1}$ when used as  preconditioners  for the pure Helmholtz matrix $A$ 
(hence  aiming to satisfy criteria (i) and (ii) of \S \ref{sec:Intro} simultaneously). The variants of \eqref{eq:defAS} which we will consider  include the Restricted, Hybrid and local impedance preconditioners. These are   
defined  in \S \ref{sec:Variants}.  

 First we give a summary of the theoretical results for \eqref{eq:defAS}. 
These are taken from \cite{GrSpVa:15}. The proofs are based on  
an  analysis of projection operators onto subspaces with 
respect to the sesquilinear form which underlies the shifted 
problem \eqref{eq:PDE}. This type of analysis is well-known for coercive 
elliptic     problems, but \cite{GrSpVa:15} was the first to devise such a theory for the high-frequency Helmholtz equation. 

%\ednote{describe technical conditions on coarse mesh and subdomains, overlap parameter $\delta$ etc}

\section{Main Theoretical Results} 
\label{sec:Main}

% Describe some of the theory from our work and its advantages and limitations. 

% The main tools for the  theory  are a  $k$- and $\eps$-explicit 
% coercivity result for  the underlying  sesquilinear form  and $k-$ and 
% $\eps-$ explicit stability theory for \eqref{eq:PDE} and its local versions on subdomains.  

Here we describe the main results from \cite{GrSpVa:15}, namely Theorems 5.6 and 5.8 in that reference. 
% This requires 
% \begin{itemize} 
% \item explaining the duality trick for right preconditioning 
% \item Giving the theorems suitably simplified for $\eps \sim k^2$ 
% \item describing briefly what happens when $\eps < k^2$.  
% \end{itemize} 

Since the systems arising from the discretisation of  \eqref{eq:PDE} are not Hermitian we need to use 
a  general purpose solver. Here we used  GMRES. Estimates of the condition number of the 
preconditioned matrix are not then enough to predict the convergence rate of GMRES. Instead 
one has to estimate  either (i)  the 
condition of the basis of  eigenvectors of the system matrix,  or  (ii) bounds on its field of values. 
Here we take the second  approach, making use of the classical theory of \cite{EiElSc:83} (see also  
 \cite{BeGoTy:06}). A brief summary of this theory is as follows. 

Consider a nonsingular linear system  
 $
\matrixC \bfx = \bfd $  in $\mathbb{C}^n$. 
Choose an initial guess $\bfx^0$ for  $\bfx$, then introduce the residual $\bfr^0 = \bfd- C \bfx^0$ and 
the usual Krylov spaces:  
$  \mathcal{K}^m(C, \bfr^0) := \mathrm{span}\{\matrixC^j \bfr^0 : j = 0, \ldots, m-1\} \ .$
Introduce a  Hermitian positive definite  matrix $D$ and the corresponding inner product on $\C^n$: 
$
\langle \bV, \bW\rangle_D := \bW^*D\bV
$, 
and let  $\Vert \cdot \Vert_D$ denote the corresponding induced norm. 
 
For $m \geq 1$, define   $\bfx^m$  to be  the unique element of $\cK^m$ satisfying  the   
 minimal residual  property: 
$$ \ \Vert \bfr^m \Vert_D := \Vert \bfd - \matrixC \bfx^m \Vert_D \ = \ \min_{\bfx \in \mathcal{K}^m(C, \br^0)} \Vert {\bfd} - {\matrixC} {\bfx} \Vert_D  , $$
When $D = I$ this is just the usual GMRES algorithm, and we 
write $\Vert \cdot \Vert = \Vert \cdot \Vert_I$, but for  more general  $D$ it 
is the weighted GMRES method \cite{Es:98} in which case  
its implementation requires the application of the weighted Arnoldi process \cite{GuPe:14}. 
The reason for including weighted GMRES in the discussion will become clear later in this section. 
% In  \S \ref{sec:Experiments} we give results for standard  GMRES which corresponds to the case $D= I$ and also for a weighted variant with respect to a certain matrix $D$ defined by \eqref{eq:weight} below.  

The following theorem  is then a simple generalisation of the classical convergence result 
stated (for $D = I$) in  \cite{BeGoTy:06}.   A proof is given in \cite{GrSpVa:15}. 

\begin{theorem} \label{Thm:A1}  Suppose    $0 \not \in W_D(\matrixC)$.  \  Then
\begin{equation}\label{eq:GMRESest}
\frac{\Vert \bfr^m \Vert_D} { \Vert \bfr^0 \Vert_D} \ \leq  \ \sin^m (\beta)   \ , 
\quad \text{where} \quad  \cos(\beta) : =
\frac{\mathrm{dist}(0,W_D(\matrixC))}{\Vert \matrixC\Vert_D }\ ,  \end{equation}
where $W_D(C)$ denotes the  {\em field of values} (also called the {\em numerical range} of $C$) with respect to the inner product induced by $D$, i.e.  
$$W_D(C) \ = \ \{\langle \bfx,\matrixC \bfx\rangle_D: \bfx \in \mathbb{C}^n, 
\Vert \bfx \Vert_D = 1\}. $$ 
\end{theorem}
This theorem shows that if the preconditioned matrix has a bounded norm,  
and has field of values bounded away from the 
origin,  then GMRES will converge independently of all parameters which are not present in the bounds.  

With this criterion  for robust convergence in mind,  the following results were  proved in \cite{GrSpVa:15}.
 These results use   the notation $A\lesssim B$ (equivalently $B \gtrsim A$) to mean  that $A/B$ is 
bounded above by a constant independent of $k$, $\eps$, and mesh diameters $h, 
\Hsub, H$.  
We write $A \sim B$ when  $A\lesssim B$ an  $B \lesssim  A$. 
In all the theoretical results below $k$ is assumed constant. 

{In Theorem \ref{thm:final1} and Corollary \ref{cor:final2} below, the matrix $D_k$ which appears is the stiffness matrix  
arising from discretising the energy inner product for the Helmholtz equation using  
the finite element basis. More precisely, the Helmholtz energy inner product and associated norm are  
defined by
\begin{equation}\label{eq:energy} 
(v, w)_{1,k} : = \int_\Omega \left(\nabla v . \nabla \overline{w} + k^2 v \overline{w} \right) \mathrm{d} x  , \quad  \quad \text{and} \quad 
\Vert v \Vert_{1,k} =  (v, w)_{1,k}^{1/2}\ .  
\end{equation}
For star-shaped Lipschitz domains,  the norm  $\Vert u\Vert_{1,k} $ of  the solution $u$ of the  
Helmholtz boundary-value problems \eqref{eq:PDE0}, \eqref{eq:ImpBC} (or alternatively  \eqref{eq:PDE} in the case of absorption)  
 can be estimated in terms of the data  $f$ and $g$ (measured in suitable norms) with a constant that is independent of $k$ and $\eps$ (provided $\eps$  grows no faster 
than $\mathcal{O}(k^2)$). This fact is the starting point (and a crucial  building 
block)  for the theory in \cite{GrSpVa:15}. If $\phi_\ell$ are the basis functions for the finite 
element space on the fine mesh,  then the matrix $D_k$ is defined by
$$ (D_k)_{\ell, m} = (\phi_\ell, \phi_m)_{1,k} \quad \text{for all } \quad \ell, m. 
$$   The matrix $D_k^{-1}$ appears as a weight in the result for right preconditioning in Theorem \ref{thm:final3}. These weights appear as artefacts of the method of analysis of the domain decomposition method 
which makes crucial use of the analysis of the Helmholtz equation in the energy norm.   Fortunately, 
in practice,  standard GMRES performs just as well as weighted GMRES (and is more efficient) - see Remark
\ref{rem:weight} below for more details.    
% Their 
% appearance implies that to really implement these theorems we should apply a weighted 
% (rather than standard) GMRES method. We did this (for the case of left preconditioning) but fortunately the results were almost the same as for standard GMRES, so all the experiments in the paper are for standard   
}

% The following are the main results from \cite{GrSpVa:15}. 

\begin{theorem}[Left preconditioning]\label{thm:final1}
\begin{align*}
(i)\quad \quad \Vert B_{\eps,AS}^{-1} A_\eps \Vert_{D_k}  & \ \lesssim  \ \left(\ksqeps\right) \quad \text{for all} \quad H, \Hsub.
\end{align*}
Furthermore, there exists a constant $C_1$ such that 
 \begin{align*}
 (ii)   \quad \quad 
 \vert \langle \bV, 
B_{\eps, AS}^{-1}A_\eps \bV \rangle_{D_k} \vert  & \ \gtrsim   \ 
%\left(1 + \frac{H}{\delta}\right)^{-1} 
\left(\frac{\eps }{k^2} \right)^{2}\ 
\Vert \bV \Vert_{D_k}^2, \quad\text{for all}\quad   \bV \in \C^n,
\end{align*}
when
\begin{equation} 
\label{eq:E20rpt}
\max\left\{ k\Hsub,\  kH 
%\left(1+ \frac{H}{\delta}\right)  
\left(\ksqeps \right)^2 \right\}  \ \leq \  
C_1 
%\left(1 + \frac{H}{\delta}\right)^{-1}  
\left(\frac{ \eps  }{k^2}\right).
\end{equation}
\end{theorem}

This result contains a lot of information. In particular, if $\eps \sim k^2$ and $kH, k\Hsub$ 
are uniformly bounded, then (weighted) left-preconditioned GMRES applied to systems with matrix $A_\eps$  
will converge in a parameter-independent way.  However when $\eps/k^2 \rightarrow 0$, 
 the bounds degrade. Nevertheless, numerical experiments in \cite{GrSpVa:15} (in the regime $H \sim \Hsub$) suggest there is some room to sharpen the theory: 
In particular, if $\eps \sim k^2$ the convergence of GMRES is parameter-independent even when $kH \rightarrow \infty$ quite quickly (that is much coarser coarse meshes than those predicted by the theory are possible). However if $\eps \sim k$ then there appears not to be  much scope to further reduce the coarse mesh diameter $H$.     
   
Combining Theorem \ref{Thm:A1} and Theorem \ref{thm:final1} we obtain:
\begin{corollary}[GMRES convergence for left preconditioning]\label{cor:final2}
Consider the weighted GMRES method where the residual is minimised in
the norm induced by $D_k$. 
Let $\br^m$ denote the $m$th iterate of GMRES applied to the system $A_\eps$, left  preconditioned with 
$B_{\eps,AS}^{-1}$. Then 
\begin{equation}\label{eq:conv_est}
\frac{\Vert \bfr^m \Vert_{D_k}} { \Vert \bfr^0 \Vert_{D_k} } \
\lesssim   \ \left(1- 
%\left(1 + \frac{H}{\delta}\right)^{-2}\, 
\left(\frac{\eps }{k^2}\right)^6\right)^{m/2} \ ,  
\end{equation}
provided condition \eqref{eq:E20rpt} holds.
\end{corollary}

Nowadays both left- and right- preconditioning play important roles in system solvers, and, 
in particular,  right preconditioning is necessary if one wants to use Flexible GMRES (FGMRES) \cite{Sa:83}. 
Fortunately Theorem \ref{thm:final1} can  be adapted to the case of right preconditioning as follows. 

The first observation is that, for any $n \times n$ complex matrix $C$ (and   
working in the inner product $\langle \cdot, \cdot \rangle_D $  induced by some SPD matrix $D$),  we have, for any $\bv \in \C^n$ and $\bw := D\bv$,   
\begin{equation} \label{eq:identity}  \frac{\langle \bv, C\bv \rangle_D}{\langle \bv, \bv \rangle_D}
\ = \   \frac{\overline{\langle \bw, C^* \bw\rangle}_{D^{-1}}}{\langle \bw, \bw \rangle_{D^{-1}}}
\ ,   \end{equation}
where $C^{\ast} = \overline{C}^\top$ denotes  the Hermitian transpose of $C$.  
Thus  estimates for the distance of the field of values of $C$ from the origin with respect to $\langle \cdot, \cdot \rangle_D$  are equivalent to  
analogous  
estimates for the field of values of $C^*$  with respect to $\langle \cdot, \cdot \rangle_{D^{-1}}$. 
   
The second   observation is  that  Theorem \ref{thm:final1} also holds for the adjoint of problem \eqref{eq:PDE}. In the adjoint case,   the sign of $\eps$ is reversed in the PDE and the boundary condition is replaced by 
$ \partial u / \partial n + \ri k u  = g$. In this case the estimates in Theorem    \ref{thm:final1} continue to hold, but  with $\eps$ replaced by $\vert \eps \vert$. This is also proved in \cite{GrSpVa:15}. 

To handle the  right-preconditioning case,  we  consider the field of values of the 
matrix $A_\eps B_{\eps, AS}^{-1}$ in the inner product induced by $D_k^{-1}$. By \eqref{eq:identity}  
these are  provided by estimates of the field of values of 
  $B_{\eps, AS}^{-*} A_{\eps}^* $ in the inner product induced by $D_k$.  The latter  are provided 
directly by the  (the extended version of)  Theorem \ref{thm:final1}. The required estimates for the norm of $A_\eps B_{\eps, AS}^{-1}$ are obtained by a similar argument. 

% Unfortunately the fact that we have to work in the weighted inner product causes 
% a slight complication. It turns out that because the estimates in Theorem \ref{thm:final1} 
% are obtained in the inner product induced by $D_k$, then the corresponding estimates for  
% the right-preconditioned matrix are obtained in the inner product induced by   $D_k^{-1}$. 
The result (from \cite{GrSpVa:15}) is as follows. 
% Since the estimates in Theorem \ref{thm:final1} are in the inner product induced by $D_k$ it follows (via the identity \eqref{eq:identity}) that the result in Theorem \ref{thm:final3} is in the inner product induced by $D_k^{-1}$.    
 
% Some of our experiments below use right preconditioning rather than
% left preconditioning. 

% Nevertheless, using the coerciveness for the adjoint
% form in Corollary \ref{cor:adjoint}, we can obtain the following  result about right
% preconditioning, however in the inner product induced by
% $D_k^{-1}$. From this, the analogue of Corollary
% \ref{cor:final2}, with $D_k$ replaced by $D^{-1}_k$, follows.

\begin{theorem}[right preconditioning] \label{thm:final3}
With the same notation as in Theorem \ref{thm:final1}, we have  
\begin{align*}
(i)\quad \quad \Vert  A_\eps B_{\eps, AS}^{-1}  \Vert_{D_k^{-1}}  & \ \lesssim  \ \left(\ksqeps\right) \quad \text{for all} \quad H, \Hsub.
\end{align*}
Furthermore, provided condition \eqref{eq:E20rpt} holds, 
 \begin{align*}
(ii) \quad\quad\vert \langle \bV, 
A_\eps B_{\eps, AS}^{-1} \bV \rangle_{D_k^{-1}} \vert  & \ \gtrsim   \ 
%\left(1 + \frac{H}{\delta}\right)^{-1} 
\left(\frac{\eps}{k^2}\right)^{2}\ 
\Vert \bV \Vert_{D_k^{-1}}^2, \quad\text{for all} \quad  \bV \in \C^n.
\end{align*}
\end{theorem}

\begin{remark} \label{rem:weight}
As described earlier, the estimates above are in the weighted inner products induced by $D_k$ and $D_k^{-1}$. It would 
be inconvenient to have to implement GMRES with these weights, especially the second one.   
It is thus  an interesting  question whether the use of weighted GMRES is necessary in practice for 
these problems. We investigated both standard and weighted GMRES (in the case of left preconditioning and with weight $D_k$) for a range of problems (some covered by the theory, some not). In practice there was little difference between the two methods.  Therefore, the   numerical experiments reported here  
use  standard GMRES. 
\end{remark} 
\begin{remark}
\label{rem:overlap}
The theorems  in \cite{GrSpVa:15} also allowed general parameter  
$\delta>0$ which described the amount of overlap between subdomains, 
and included the dependence on $\delta$ explicitly in the estimates. We suppressed this here in order to make the 
exposition simpler.    
\end{remark}     
% The DD preconditioner studied theoretically in \cite{GrSpVa:15} 
% is the conventional Schwarz method (with Dirichlet conditions on subdomain boundaries). This is always well-defined for $\eps>0$ because of a coercivity 
% property for the underlying operator (see \cite[Lemma 2.4]{GrSpVa:15}). 
% Let's denote the question  of performance of  preconditioners for \eqref{eq:PDE} with $\eps>0$ as ``Problem A''. 

% The  estimates in \cite{GrSpVa:15}  are obtained with respect to  
%  the Helmholtz energy norm:
% \begin{equation} \label{eq:HE}  \Vert v \Vert_{1,k}^2 : = \vert v \vert_{H^1(\Omega)}^2 + \omega^2 \Vert v \Vert_{L^2(\Omega)}^2 \ ,  
% \end{equation} 
% and the field of values is defined with respect to the corresponding energy inner product. That means there is a gap between the proofs and the  
% standard Euclidean norm implementation of GMRES used in normal practice 
% (a phenomenon also studied in the elliptic-dominated case in \cite{CaZo:02}).   
\section{Variants of the Preconditioners} 
\label{sec:Variants}

In this section we describe the variants of the classical Additive Schwarz method defined in \eqref{eq:defAS} which are investigated in the numerical {experiments} which follow.  
% Each of these  variants  replace the
% local component  \eqref{eq:minor} with something else. 
% The first variant involves averaging   
%  in the overlap of the subdomains. For each fine grid node $x_j  \ \ (j \in \cI^h)$,  let
%  $L_j$ denote the number of subdomains which contain  $x_j$. 
% Then the local operator is:  
% $$
% (\matrixB_{\eps,AVE,local}^{-1}\bfv)_j  \ = \  \frac{1}{L_j} \sum _{\ell: \, x_j \in {\Omega_\ell}}
% \left(R_\ell^T A_{\eps,\ell}^{-1} R_\ell\bfv\right)_j\ , \quad
% \text{for each} \quad j \in \cI^h\ , 
% $$
% and  the corresponding 
%   {\em Averaged Additive Schwarz} ({\tt AVE}) preconditioner   is: 
% \begin{equation}\label{eq:defAVE}
% \matrixBepsAVE  \ = \ R_0^T A_{\eps,0}^{-1}
% R_0  \  +\  \matrixB_{\eps,AVE,local}^{-1} . 
% \end{equation}

The first  variant which we consider  is  the {\em Restrictive Additive Schwarz} ({RAS})
 preconditioner, which is well-known in the literature  \cite{CaSa:99}, \cite{KiSa:13}. Here,   
to define the local operator,  
 for each $j \in \cI^h$, choose a single $\ell =
\ell(j)$ with the property that $x_j \in \Omega^{\ell(j)}$.  
Then  
the action of the local contribution,  for each vector of fine grid freedoms  $\bfv$, is: 
  \begin{equation}
\label{eq:RASlocal}
(B_{\eps,RAS,local}^{-1}\bfv)_j  \ = \ 
\left(R_{\ell(j)}^T A_{\eps,\ell(j)}^{-1} R_{\ell(j)}\bfv\right)_j\ , \quad
\text{for each} \quad j \in \cI^h\ .
\end{equation}
{We denote this one level preconditioner as RAS1. (We shall in fact use a slight variation on this -  as  described precisely in  \S \ref{sec:Experiments}.)   }  

From this we could build the RAS preconditioner (in analogy to the standard Additive Schwarz method):   
\begin{equation}\label{eq:defRAS}
\matrixBepsRAS  \ = \ R_0^T A_{\eps,0}^{-1}
R_0  \  +\  \matrixB_{\eps,RAS,local}^{-1}\ .
\end{equation}

However we shall not use this directly in the following. Rather,  instead 
of doing all the local and 
coarse grid problems independently 
(and thus potentially in 
parallel), we  first do  a 
coarse solve and then perform the   local solves on the residual of the coarse solve. 
This was first introduced in  \cite{MaBr:96}. As described in \cite{GrSc:07},   
this method is
closely related to the deflation method \cite{NaVu:04},  which has been used recently
to good effect in the context of shifted Laplacian combined with multigrid 
\cite{ShLaVu:13}. 
The Hybrid RAS ({HRAS}) preconditioner then takes  the form   
\begin{equation}\label{eq:HRAS}
\matrixB_{\eps,{HRAS}}^{-1} :=  \matrixR_0^T \matrixA_{\eps,0}^{-1} \matrixR_0  +  \matrixP_0^T \left(  \matrixB_{\eps,RAS,local}^{-1} \right)  \matrixP_0 
\  , \end{equation}
where 
 $$\matrixP_0 \ = \ \matrixI - \matrixA  \matrixR_0^T
 \matrixA_{\eps,0}^{-1} \matrixR_0\ . $$

Remembering that the local solves in $B_{\eps,RAS,local}^{-1}$ are
solutions of local problems with a Dirichlet condition on interior
boundaries of subdomains, and noting that these are not expected to perform 
well for genuine wave propagation (i.e.~$\eps$  small and $k$ large), 
we also  consider the use of impedance boundary  conditions on the local solves.
% We therefore  
% introduce the  sesquilinear form local to the subdomain $\Omega_\ell$,
% defined as the following local equivalent of
% \eqref{eq:Helmholtzvf_intro}:  
% $$ a_{\eps,Imp,\ell}(v,w) \ = \ \int_{\Omega_\ell} \nabla v . \nabla
% \overline{w}  -    (k^2 + i \eps) \int_{\Omega_\ell} v \overline{w} - \ri
% k \int_{\partial \Omega_\ell} v \overline{w}, 
% $$
% (remember that we are choosing $\eta = k$ in all the experiments). We
Let $A_{\eps,Imp,\ell}$ be the stiffness matrix arising from the solution of \eqref{eq:PDE} restricted to  
$\Omega_\ell$, where the impedance condition $\partial u / \partial n - \ri k u$ is imposed on the boundary $\partial \Omega_\ell$, and dealt with in the finite element method as a natural boundary condition.  
%\ednote{gave more detail on this} 
% $$(A_{\eps,Imp,\ell})_{j,j'} \ = \ a_{\eps,Imp, \ell}(\phi_{j'}, \phi_{j})\ ,
% \quad j,j' \in \cI(\overline{\Omega_\ell}) .$$ 
This can be used as a local operator in the  {\tt HRAS} operator 
\eqref{eq:HRAS}. The one-level variant is  
  \begin{equation}
\label{eq:ImpRASlocal}
(\matrixB_{\eps,Imp,RAS,local}^{-1}\bfv)_j  \ = \ 
\left(\widetilde{R}_{\ell(j)}^T A_{\eps,Imp,\ell(j)}^{-1} \widetilde{R}_{\ell(j)}\bfv\right)_j\ , \quad
\text{for each} \quad j \in \cI^h\ ,
\end{equation}
{Here (noting that the local impedance condition is handled as a natural  boundary condition on $\Omega_\ell$), $\widetilde{R}_\ell$ denotes the 
restriction operator $(\tilde{R}_{\ell})_{j,j'} = \delta_{j,j'}$,  (as before)   $j'$ ranges  
over all  $\cI^h$, but now   
$j$ runs over all  indices such that  $x_j \in \overline{\Omega}_\ell$. } 

The hybrid two-level variant is 
\begin{equation}\label{eq:HRASImp}
\matrixB_{\eps,{Imp, HRAS}}^{-1} :=  \matrixR_0^T \matrixA_{\eps,0}^{-1} \matrixR_0  +  \matrixP_0^T \left(  \matrixB_{\eps,Imp,RAS,local}^{-1} \right)  \matrixP_0 
\  . \end{equation}
We refer to these as the one- and two-level {ImpHRAS} preconditioners. 
%\vspace{0.5cm} 

In the following section we will concentrate on illustrating the use of the four  preconditioners 
defined in \eqref{eq:RASlocal}, \eqref{eq:HRAS}, 
\eqref{eq:ImpRASlocal}  and \eqref{eq:HRASImp} for solving various  problems with system matrix $A$ (i.e. the discretisation of \eqref{eq:PDE} with $\eps = 0$). 
In our discussion and in the  tables below  we will use  the following notation for the preconditioners: 
\begin{equation}\label{notation} 
\eqref{eq:RASlocal}  = \text{RAS1} , \quad \quad \eqref{eq:HRAS} = \text{HRAS}, \quad \quad 
\eqref{eq:ImpRASlocal} = \text{ImpRAS1},  \quad \quad  \eqref{eq:HRASImp} = \text{ImpHRAS} \ .
\end{equation}

% denote these four preconditioners by RAS1, HRAS, ImpRAS1 and ImpHRAS. 

% $$ B_{\eps,RAS,local}^{-1}, \quad B_{\eps,HRAS}^{-1}, \quad    B_{\eps,Imp,RAS,local}^{-1}, \quad  
% B_{\eps,{Imp, HRAS}}^{-1}$$

\section{Numerical Experiments} 
\label{sec:Experiments} 

% This section describes the benefits (and some limitations) in using local impedance problems in DD methods for the Helmholtz equation. 

Our numerical experiments concern the solution of \eqref{eq:PDE} on
the unit square, with $\eta = k$ and $\eps = 0$,  
discretised by the continuous piecewise linear finite element
method on a uniform triangular mesh.
Thus,  the problem being solved here is the ``pure Helmholtz''  problem without absorption and can be completely  specified   by  the  {\em fine mesh diameter}, here denoted  
   $\hprob$. In \cite{GrSpVa:15} we also computed iteration numbers for solving 
\eqref{eq:PDE} with  $\eps> 0$, thus  an additional parameter $\epsprob$ was   needed 
to specify the problem being solved.  
% The case $\epsprob>0 $ was studied in \cite{GrSpVa:15} in order to verify  the 
%  predictions of  Theorems \ref{thm:final1} and \ref{thm:final3}. 
Here  we restrict to the case  $\epsprob = 0$. 
For the solver  we shall use domain decomposition preconditioners built from 
various approximate inverses for \eqref{eq:PDE}. The choice of $\eps>0$ 
which is used to build the preconditioner is denoted $\epsprec$.   

%In the case of constant wavenumber $k$, 
% We will discuss both the 
% cases $\eps \not = 0$ and $\eps = 0$.   

The experiments in \S \ref{subsec:pollution-free}  will be concerned with the case    when the fine grid diameter is   
$\hprob \sim   k^{-3/2}$. This is  
the discretisation level  generally believed to be necessary to  remove the pollution effect: roughly speaking 
the relative error obtained with this choice of $\hprob$ is not expected  to grow as $k \rightarrow \infty$. (However there 
is no proof of this except in the 1D case: See, e.g., the literature reviews in \cite[Remark 4.2]{GaGrSp:15} and
\cite[\S1.2.2]{GrLoMeSp:15}.) 

 However the case of a fixed number of grid points per wavelength ($\hprob \sim k^{-1}$) is also frequently used 
in practice (especially in 3D) and  provides  sufficient accuracy in a limited  frequency range. 
This regime is often studied in  papers about 
Helmholtz solvers and so we include a substantial subsection (\S \ref{subsec:pollution}) 
on results for  this case, which was  not specifically discussed  in \cite{GrSpVa:15}.  
Nevertheless the question  of preconditioning the problem defined by  
$\hprob \sim k^{-1}$ and $\epsprob \sim k$ did arise in \cite{GrSpVa:15}, 
as an ``inner problem'' in the multilevel solution of the problem 
with $\hprob \sim k^{-3/2}$,  $\epsprob = 0$. (This is discussed again 
in \S \ref{subsec:pollution-free} below.)  
  
Interestingly, it turns out that the asymptotics (as $k$ increases)  of the solvers in each of 
the two cases $\hprob \sim k^{-3/2}$ and $\hprob \sim k^{-1}$ (both with $\epsprob = 0$) 
are somewhat 
different from each other and the best methods for one case  are not necessarily the best for the other.

In the general theory given in \S \ref{sec:Main},  
coarse grid size $H$ and subdomain  size $\Hsub$  
are  permitted to be unrelated. In our  
experiments here we construct  
local subdomains  
by first choosing a coarse grid and then taking each of the elements of the 
coarse grid and extending them to obtain an
overlapping cover of subdomains with overlap parameter $\delta$. This is chosen  
as large as possible, but with the restriction no two extended subdomains can touch unless they came from touching elements of the original coarse grid.   In the literature this is called 
{\em generous overlap} and $\Hsub \sim H$.  
Thus our preconditioners are completely  determined  by specifying the values  
of  $H$ and $\eps$. In the case of constant $k$,  we denote these  by 
\begin{equation}\label{eq:precspec}
\Hprec  \quad \text{and} \quad \epsprec\ . 
\end{equation}

We also have to specify how the RAS subdomains (recall \eqref{eq:RASlocal}) are defined.
Actually in our implementation involves a slight variation on    \eqref{eq:RASlocal} as follows. 
Our  RAS subdomains are the original elements of the coarse grid (before extension). 
These overlap, but  only at the edges of the coarse grid. 
Each node of the fine grid lies in a unique  
RAS  subdomain except for nodes on the coarse grid edges. 
At these nodes the RAS operator \eqref{eq:RASlocal} is extended so that it 
performs averaging of the contributions from all relevant subdomains at all such edge nodes.

When designing good domain decomposition methods we should be aware of cost. In the 
classical context (which we adopt here)  where coarse grid and local problems are 
linked, a large-sized coarse grid problem  will imply small-sized local problems and vice-versa. Coarse grids which are very fine and very coarse can both lead to very good methods in terms of iteration numbers, but not necessarily optimal in terms of time. 

An ``ideal'' situation may be when all sub-problems  are ``load balanced''.
Let $\hprob$ be the  fine grid diameter  and let $\Hprec$ be the coarse grid diameter, 
so that in $\mathbb{R}^d$, the dimension of the coarse grid problem is $\mathcal{O}(\Hprec^{-d})$, while   
the dimension of the local problems are $\mathcal{O}((\Hprec/\hprob)^d)$. 
Then  the classical domain decomposition  method  is load-balanced when $\Hprec \sim  \hprob^{1/2}$. 
If generous overlap is used, then  a slightly  smaller $\Hprec$ will 
give us load balancing.   For example, in the pollution-free case    
$\hprob = k^{-3/2}$, the domain decomposition  will be load-balanced at 
about $\Hprec = k^{-0.8}$.  While load balancing occurs at about  $\Hprec\sim k^{-0.6}$ when we are taking a fixed number of points per wavelength  ($\hprob \sim k^{-1}$).  We use these estimates as a guide in the experiments below.  

% In the reminder of this section we perform experiments on problem \eqref{eq:PDE} 
%and \eqref{eq:ImpBC} in the case when $\Omega$ is the unit square $[0,1]\times [0,1]$.     

In all the experiments below the stopping tolerance for GMRES was that the 
relative residual should be reduced by $10^{-6}$. 

{In the experiments below, the system being solved is always the pure Helmholtz system $A \mathbf{u} = \mathbf{f}$. In the results given in Tables 1-3 the right hand side vector $\mathbf{f}$ was chosen so that the finite element solution is an approximation of a plane wave (see \cite[\S6.2]{GrSpVa:15}). For the rest of the experiments $\mathbf{f} = \mathbf{1}$ was used.    }

\subsection{Pollution-free systems  ($\hprob \sim k^{-3/2}$)} 
\label{subsec:pollution-free}
{The timings given in Tables \ref{tab:1} - \ref{tab:3}  below  were for implementation  
on a serial workstation with Intel Xeon E5-2630L CPUs with 48GB RAM. The later experiments were on a multiprocessor, described  in \S 
\ref{subsec:pollution}. } 
% The additive domain decomposition algorithms described below are highly
% parallelisable and parallel scaling will be investigated in a later paper.
% \ednote{Change this if Eero's results work out {\color{blue} Note that Table 4
% and 5 were run on a different cluster. (Details in Section 5.2.1.)}} 

The performance of GMRES for this case is investigated in detail in 
\cite{GrSpVa:15}. There we first studied the performance of domain 
decomposition  preconditioners for systems with absorption 
(i.e. we set $\epsprob = \eps >0$ and we studied  the performance of 
$B_\eps^{-1}$ as a preconditioner for $A_\eps$). With respect to that question we found  that: \\

\noindent 
(i) the performance of the solvers  reflected the theory given in \S \ref{sec:Main}; \\
(ii) There was little difference between left- and right-preconditioning;\\
(iii) There was little difference between the performance of standard GMRES and GMRES which minimised the residual in the weighted norm (in the case of left preconditioning) induced by $D_k$  (see Remark \ref{rem:weight} at the end of \S \ref{sec:Main});\\   
(iv) There was 
a marked superiority for HRAS over several other   variants of Additive Schwarz; \\
(v) If $\Hprec$ is small enough ($\Hprec\sim k^{-1}$ is sufficient),  then $B_\eps^{-1}$ is a good  
preconditioner for $A_\eps$ even for rather small $\eps$ (in fact,  even $\eps = 1$ gives acceptable results for HRAS);\\ 
(vi) If $\Hprec$ is small enough then it makes little difference whether the local problems have Dirichlet or impedance boundary conditions;  
\\
(vii) For  larger $\Hprec$, Dirichlet local problems perform very badly, while impedance local problems  work well for large enough $\Hprec$. In this case the coarse grid solver can be switched off without degrading the convergence of GMRES. \\

Based on these observations, the discussion in \cite{GrSpVa:15} then turned 
to the more important question of the solution of problems without absorption (i.e. $\epsprob = 0$).
The  discussion in the rest of this subsection  is an expansion of the discussion in \cite{GrSpVa:15}.        

 We compare HRAS (Hybrid Restricted Additive Schwarz  
with Dirichlet local problems), as defined in \eqref{eq:HRAS}  with  ImpHRAS 
(Hybrid RAS with Impedance local problems), as defined in \eqref{eq:HRASImp}.
In these experiments, $\hprob = k^{-3/2}$ and in Table \ref{tab:1} below  
we give the number of GMRES iterations (with \# denoting  iteration count) for each of these two methods 
for various choices of $\Hprec$ and $\epsprec$.   In Table \ref{tab:1},  the headline figure for each case is the iteration number for the Hybrid method \eqref{eq:HRAS}
or \eqref{eq:HRASImp}, while as a subscript we give the  iteration count  
for the corresponding  one level methods (omitting the coarse grid solve), given respectively by \eqref{eq:RASlocal} and \eqref{eq:ImpRASlocal}. We include iteration numbers for the three 
cases $\epsprec = k, k^{1.2}, k^2$. The   optimal choice turns out to be  
around $\epsprec \in  [k, k^{1.2}]$, 
while $\epsprec = k^2$  is provided for comparison.  Data for a larger range of $\epsprec$ and $\Hprec$ is given in \cite{GrSpVa:15}. 
  A $*$ in the tables means the iteration did not converge after 200 iterations.  
%\ednote{Eero is this correct?}    
\begin{table} 
\begin{center} 
\begin{tabular}{ccc}
\begin{tabular}{|c|c|c|}
\hline 
\multicolumn{3}{|c|}{\quad $\Hprec \sim k^{-1}$, \ $\epsprec = k$ \quad } \tabularnewline
\hline 
$k$  & \# HRAS & \# ImpHRAS \tabularnewline
\hline 
20 & $12_{92}$ & $17_{105}$ \\
40 & $18_{*}$ & $21_{*}$\\
60 & $25_{*}$ & $27_{*}$\\
80& $33_{*}$ & $35_{*}$ \\
100 & $43_{*}$ & $45_{*}$\\
\hline
\end{tabular} 
&
\begin{tabular}{|c||c|c|}
\hline 
\multicolumn{3}{|c|}{\quad $\Hprec \sim k^{-1}$, \ $\epsprec = k^{1.2}$ \quad } \tabularnewline
\hline 
$k$  & \# HRAS & \# ImpHRAS \tabularnewline
\hline 
20 & $13_{92}$ & $18_{105}$ \\
40 & $18_{*}$ & $21_{*}$\\
60 & $25_{*}$ & $27_{*}$\\
80& $32_{*}$ & $34_{*}$ \\
100 & $42_{*}$ & $43_{*}$\\
\hline
\end{tabular} 
&
\begin{tabular}{|c||c|c|}
\hline 
\multicolumn{3}{|c|}{\quad $\Hprec \sim k^{-1}$,\ $\epsprec = k^2$ \quad } \tabularnewline
\hline 
$k$  & \# HRAS & \# ImpHRAS \tabularnewline
\hline 
20 & $37_{93}$ & $34_{113}$ \\
40 & $63_{*}$ & $56_{*}$\\
60 & $86_{*}$ & $78_{*}$\\
80& $110_{*}$ & $101_{*}$ \\
100 &$136_{*}$  & $123_{*}$\\
\hline
\end{tabular} 
\end{tabular} 
% \caption{Comparison of HRAS and ImpHRAS when  $\hprob \sim k^{-3/2}$,  
% $\epsprob = 0$, and $\Hprec \sim k^{-1}$ } 
\end{center} 
%\end{table}

%\begin{table} 
\begin{center} 
\begin{tabular}{ccc}
\begin{tabular}{|c|c|c|}
\hline 
\multicolumn{3}{|c|}{ \quad $\Hprec \sim k^{-0.6}$,\ $\epsprec = k$ \quad } \tabularnewline
\hline 
$k$  & \# HRAS & \# ImpHRAS \tabularnewline
\hline 
20 & $51_{63}$ & $26_{31}$ \\
40 & $125_{133}$ & $50_{51}$\\
60 & $*_{*}$ & $69_{71}$\\
80& $*_{*}$ & $74_{84}$ \\
100 & $*_{*}$ & $84_{97}$\\
\hline
\end{tabular} 
&
\begin{tabular}{|c||c|c|}
\hline 
\multicolumn{3}{|c|}{ \quad $\Hprec \sim k^{-0.6}$,\ $\epsprec = k^{1.2}$ \quad } \tabularnewline
\hline 
$k$  & \# HRAS & \# ImpHRAS \tabularnewline
\hline 
20 & $48_{58}$ & $26_{32}$ \\
40 & $114_{125}$ & $48_{51}$\\
60 & $*_{*}$ & $69_{70}$\\
80& $*_{*}$ & $74_{83}$ \\
100 & $*_{*}$ & $84_{95}$\\
\hline
\end{tabular} 
&
\begin{tabular}{|c||c|c|}
\hline 
\multicolumn{3}{|c|}{ \quad $\Hprec \sim k^{-0.6}$,\ $\epsprec = k^2$ \quad } \tabularnewline
\hline 
$k$  & \# HRAS & \# ImpHRAS \tabularnewline
\hline 
20 & $39_{43}$ & $36_{42}$ \\
40 & $81_{6}$ & $73_{66}$\\
60 & $113_{102}$ & $104_{91}$\\
80& $135_{121}$ & $126_{111}$ \\
100 & $156_{141}$ & $148_{131}$\\
\hline
\end{tabular} 
\end{tabular} 
\caption{Comparison of HRAS and ImpHRAS for the problem with   $\hprob \sim k^{-3/2}$,  
$\epsprob = 0$, using  various choices of $\Hprec$ and $\epsprec$,  \label{tab:1}}
\end{center} 
\end{table}

% \begin{table} 
% \begin{center} 
% \begin{tabular}{ccc}
% \begin{tabular}{|c|c|c|}
% \hline 
% \multicolumn{3}{|c|}{$H = k^{-1}$} \tabularnewline
% \hline 
% $k$  & HRAS & ImpHRAS \tabularnewline
% \hline 
% 40 & $18_{*}$ & $21_{*}$\\
% 60 & $25_{*}$ & $27_{*}$\\
% 80& $32_{*}$ & $34_{*}$ \\
% 100 & $42_{*}$ & $43_{*}$\\
% \hline
% \end{tabular} 
% &
% \begin{tabular}{|c||c|c|}
% \hline 
% \multicolumn{3}{|c|}{$H = k^{-0.8}$} \tabularnewline
% \hline 
% $k$  & HRAS & ImpHRAS \tabularnewline
% \hline 
% 40 & $96_{*}$ & $49_{106}$\\
% 60 & $169_{*}$ & $99_{150}$\\
% 80& $*_{*}$ & $168_{193}$ \\
% 100 & $*_{*}$ & $*_{*}$\\
% \hline
% \end{tabular} 
% &
% \begin{tabular}{|c||c|c|}
% \hline 
% \multicolumn{3}{|c|}{$H = k^{-0.6}$} \tabularnewline
% \hline 
% $k$  & HRAS & ImpHRAS \tabularnewline
% \hline 
% 40 & $48_{58}$ & $48_{51}$\\
% 60 & $114_{125}$ & $69_{70}$\\
% 80& $*_{*}$ & $74_{83}$ \\
% 100 & ?  & $84_{95}$\\
% \hline
% \end{tabular} 
% \end{tabular} 
% \caption{Comparison of HRAS and ImpHRAS with when  $\epsprob = 0$, $\epsprec = k^{1.2}$ and  
% $h \sim k^{-3/2}$ } 
% \end{center} 
% \end{table}

%\vspace{0.5cm} 
\noindent 
Based on the results in Table \ref{tab:1},   we can  make the following observations: \\

\noindent 
(i) When $\Hprec \sim  k^{-1}$,  the coarse grid is sufficiently fine and does a good job. 
Using the data for $\Hprec \sim k^{-1}$ and $\epsprec \sim k^{1.2}$ we observe that  we have 
$\# HRAS  \sim k^{0.71}$. Since we are here solving problems of size $n \sim k^3$, this is equivalent to $\# HRAS  \sim n^{0.24}$.    (Throughout the paper,   
rates of growth are obtained by linear least squares fits to the relevant log-log data.) 
Note that when $\Hprec \sim  k^{-1}$,  there is little difference between HRAS and ImpHRAS, i.e.   
it   does not matter here  whether 
the local problems have Dirichlet or Impedance condition.   
This preconditioner has a competitive performance as $n$ increases, but it incorporates  
an expensive coarse grid solve of size $~\Hprec^{-2} \sim  k^2$ and it does not work without the coarse solve.  
  \\    
(ii)  When $\Hprec \sim  k^{-0.6}$ the local problems are rather large 
(size  $ \sim  k^{9/5}$) the ImpHRAS method works reasonably  well with a 
slow growth of iteration count with respect to $k$ (although higher actual iterations), while HRAS is not usable. Moreover  in the case of 
ImpHRAS, the coarse grid solve has almost no effect and can be neglected.     \\
(iii) In all cases the best choice of absorption parameter $\epsprec$  seems to be 
about $\epsprec \sim k^{\beta}$ with $\beta$ close to $1.2$. We note that the choice  
$\epsprec \sim k^2$ is  remarkably  inferior. 
A more extensive study of the variation of iteration numbers with respect to  $\epsprec$ and $\Hprec$ is given in \cite{GrSpVa:15}. 

These observations  led to the formulation of an  
inner-outer strategy for problems with $\hprob \sim k^{-3/2}$,  
with the outer iteration having preconditioner specified by     
$\Hprec =k^{-1}$ and $\epsprec = k^{1.2}$. 
% In this case 
%  the relatively large coarse grid problem of mesh diameter  $\sim k{-1}$ and with absorption 
% $\sim k$. 
This ``outer preconditioner''  is a discretisation of \eqref{eq:PDE}  
with $\hprob  \sim k^{-1}$ and $\epsprob \sim k^{1.2}$, which is to be 
solved by a preconditioned inner iteration.  So,  as a precursor to formulating the inner-outer method,  we  
study iteration 
counts for typical instances of this inner iteration. 
Here are some sample results  with $\hprob = \pi/5k \sim k^{-1} $, $\epsprob = k^{1.2}$ using ImpHRAS 
as a preconditioner,   with $\Hprec \sim  k^{-1/2}$ and $\epsprec = k^{1.2}$:
  
\begin{table} 
\begin{center}
\begin{tabular}{|c||c|}
\hline 
$k$  &  \#ImpHRAS \tabularnewline
\hline 
20&  $ 14_{16}$\\
40 &  $21_{23}$\\
60 &  $28_{30}$\\
80&  $32_{31}$ \\
100 &  $36_{34}$\\
120 & $ 39_{38}$ \\
140 & $ 43_{41}$ \\
\hline
\end{tabular}
\caption{\label{tab:2} Iteration numbers for  ImpHRAS with 
$\epsprob = k^{1.2}  = \epsprec$,  $\hprob = \pi/5k$ and $\Hprec \sim k^{-1/2}$} 
\end{center}
\end{table}

We see from  Table \ref{tab:2} that,  even without the coarse solve,  the iteration numbers 
grow slowly, and even seem to be slowing down as $k$ increases. Extrapolation using the last 
five entries of Table \ref{tab:2} (without the coarse solve)  indicates that  
\#ImpHRAS grows with approximately 
$\mathcal{O}(k^{0.38}) = \mathcal{O}(n^{0.19})$, where $n$ is the size of  
the systems being solved in Table \ref{tab:2}. 

Therefore in \cite{GrSpVa:15} we proposed an  inner-outer FGMRES  iteration  
using (as the outer solver) HRAS with $\Hprec = k^{-1}$  and (as the inner solver) ImpRAS1
 with $\Hprec = k^{-1/2}$.  
    This method solves a system of dimension $\mathcal{O}(k^3)$ by  solving 
$\mathcal{O}(k^2 + k)$  independent subdomain problems of dimension      
$\mathcal{O}({k^{1/2} \times k^{1/2}}) = \mathcal{O}(k)$ and    was found  to have competitive properties.  

In particular the subproblems are sufficiently small as to be very efficiently solved by a 
sparse direct solver.  (Here we use {\tt umfpack} included in the scipy sparse matrix package.) In this regard, an interesting 
observation is that, while positive definite systems coming 
from 2D finite element approximations of elliptic problems are often reported to be 
solvable by sparse direct solvers in optimal time 
($\mathcal{O}(n)$,  for dimension $n$  up to about $10^5$), 
this appears not 
to be the case for the indefinite systems encountered here.  In our experience 
the computation time for the sub-systems encountered here grows  slightly faster than  linearly 
with respect to  
 dimension $n$.

% but the number of iterations still grows with $k$, 
% although the local solves needed in the implementation of the preconditioner 
% are all quite small (of size $k$).   
% The study of the composite multilevel algorithm in \cite{GrSpVa:15} 
% proposed to use $\eps = k$ in  
% the preconditioner. 
%Further  multilevel extensions are easy to envisage.   

The following table gives some sample results for the composite 
inner/outer algorithm with $\epsprec = k^{\beta}$ {(for both inner and outer iterations,  for various $\beta$)}  and 
an inner tolerance $\tau = 0.5$ (found in \cite{GrSpVa:15} to be empirically 
best).  
The numbers in bold font denote the  number of outer (respectively inner) iterations, 
while the 
smaller font numbers underneath denote the  total time in seconds 
[with an average time for each outer iteration in square brackets]. 
{(Other choices of inner tolerance are explored in  \cite{GrSpVa:15}. Recall that the outer tolerance is $10^{-6}$.)}  
\begin{table}
\begin{center}
\begin{tabular}{|c||c|c|c|c|c|c|c|}
\hline 
$k\backslash\beta$  & \textbf{0 } & \textbf{0.4 } & \textbf{0.8 } & \textbf{1 } & \textbf{1.2 } & \textbf{1.6 } & \textbf{2.0}\tabularnewline
\hline 
% \multirow{2}{*}{\textbf{10} } & \textbf{17(2) } & \textbf{17(1) } & \textbf{18(1) } & \textbf{18(1) } & \textbf{18(1) } & \textbf{19(1) } & \textbf{23(1)}\tabularnewline
%  & {\scriptsize{0.61~{[}0.02{]}}} & {\scriptsize{ 0.59~{[}0.02{]}}} & {\scriptsize{ 0.66~{[}0.01{]}}} & {\scriptsize{ 0.66~{[}0.02{]}}} & {\scriptsize{ 0.65~{[}0.02{]}}} & {\scriptsize{ 0.66~{[}0.01{]}}} & {\scriptsize{ 0.65~{[}0.01{]}}}\tabularnewline
\multirow{2}{*}{\textbf{20 }} & \textbf{19(2) } & \textbf{19(2) } & \textbf{19(2) } & \textbf{19(2) } & \textbf{19(2) } & \textbf{25(1) } & \textbf{36(1)}\tabularnewline
 & {\scriptsize{3.86~{[}0.08{]}}} & {\scriptsize{   3.72~{[}0.08{]}}} & {\scriptsize{3.72~{[}0.08{]}}} & {\scriptsize{3.68~{[}0.08{]}}} & {\scriptsize{   3.66~{[}0.08{]}}} & {\scriptsize{   4.00~{[}0.07{]}}} & {\scriptsize{   4.96~{[}0.07{]}}}\tabularnewline
\hline 
\multirow{2}{*}{\textbf{40 }} & \textbf{22(4) } & \textbf{22(4) } & \textbf{22(4) } & \textbf{22(3) } & \textbf{22(3) } & \textbf{28(2) } & \textbf{61(1) }\tabularnewline
 & {\scriptsize{54.8~{[}0.73{]}}} & {\scriptsize{54.9~{[}0.73{]}}} & {\scriptsize{  54.8~{[}0.72{]}}} & {\scriptsize{  54.7~{[}0.71{]}}} & {\scriptsize{  54.8~{[}0.71{]}}} & {\scriptsize{  58.0~{[}0.69{]}}} & {\scriptsize{  80.4~{[}0.68{]}}}\tabularnewline
\hline 
\multirow{2}{*}{\textbf{60 }} & \textbf{28(5)} & \textbf{ 28(5)} & \textbf{ 28(5)} & \textbf{ 28(5)} & \textbf{ 28(4)} & \textbf{ 35(2)} & \textbf{ 82(1)}\tabularnewline
 & {\scriptsize{370~{[}3.20{]}}} & {\scriptsize{ 371~{[}3.20{]}}} & {\scriptsize{372~{[}3.19{]}}} & {\scriptsize{ 370~{[}3.16{]}}} & {\scriptsize{ 369~{[}3.11{]}}} & {\scriptsize{ 383~{[}3.00{]}}} & {\scriptsize{539 {[}3.10{]}}}\tabularnewline
\hline 
\multirow{2}{*}{\textbf{80 }} & \textbf{36(6)} & \textbf{ 36(6)} & \textbf{ 36(6)} & \textbf{ 36(5)} & \textbf{ 35(5)} & \textbf{ 42(2)} & \textbf{104(1)}\tabularnewline
 & {\scriptsize{1288~{[}8.62{]}}} & {\scriptsize{1375~{[}8.69{]}}} & {\scriptsize{1300~{[}8.59{]}}} & {\scriptsize{1316~{[}8.51{]}}} & {\scriptsize{1273~{[}8.38{]}}} & {\scriptsize{1323~{[}8.08{]}}} & {\scriptsize{1909~{[}8.19{]}}}\tabularnewline
\hline 
\multirow{2}{*}{\textbf{100 }} & \textbf{46(8)} & \textbf{ 46(8)} & \textbf{ 46(7)} & \textbf{ 45(7)} & \textbf{ 44(6)} & \textbf{ 49(2)} & \textbf{126(1)}\tabularnewline
 & {\scriptsize{3533~{[}16.5{]}}} & {\scriptsize{3678~{[}16.01{]}}} & {\scriptsize{3586~{[}16.4{]}}} & {\scriptsize{3471~{[}15.9{]}}} & {\scriptsize{3483~{[}16.2{]}}} & {\scriptsize{3503~{[}15.5{]}}} & {\scriptsize{4832~{[}16.4{]}}}\tabularnewline
\hline 
\end{tabular}
%\end{tiny}
%\end{center} 
\caption{\label{tab:3}
 GMRES iteration counts and timings for the inner-outer algorithm with $\epsprob = 0$, $\hprob = k^{-3/2}$, $\Hprec = k^{-1}$ in the outer iteration,  $\Hprec = k^{-1/2}$ in the inner iteration and  $\epsprec = k^{\beta}$ in both inner and outer iterations} 
\end{center} 
\end{table}
The best results  occur with   $\epsprec = k^{\beta}$ with  $\beta \in [1,1.2]$. Using 
the data in the column headed $\beta = 1$ (and remembering that we are here 
solving systems of dimension $n = k^3$),     the outer  iteration count grows with about  
$\mathcal{O}(k^{0.53}) \approx \mathcal{O}(n^{0.18})$, while the time per iteration is  about  $ \mathcal{O}(n^{1.11})$ and the total time is $ \mathcal{O}(n^{1.43})$.   
 {To give an idea of the size of the systems being solved, when  $k = 100$, $n = 1,002,001$. }

An interesting observation in Table \ref{tab:3} is 
the relative insensitivity of the results to the choice of $\beta$ in the range 
$\beta \in [0,1.6]$, and the very poor performance of $\beta = 2$. Thus for this method the 
choice of absorption $\epsprec = k^2$ is a relatively poor one, while in fact  the choice $\epsprec = 1 = k^0$  is quite  competitive. This is quite different to the experience reported using multigrid shifted Laplacian preconditioners. Note also that the number of inner iterations decreases as we read the rows of Table \ref{tab:3} from left to right, because increasing $\beta$ means putting more absorption into the preconditioner and hence makes the inner problem easier to solve.  

{The remainder of the experiments in the paper 
were done on a linux cluster of 130 nodes. Each 
node consists of 2 CPUs (Intel Xeon E5-2660 v2 @ 2.20GHz) with 10 cores: in total 20 cores and  
64GB RAM on each node. The nodes are  connected with 4x QDR Infiniband networks. 
This cluster was used in serial mode except for the modest  parallel experiment in Table \ref{tab:4p}, in which up to 10 of the 130 nodes were used.  } 

\subsection{$10$ grid-points per wavelength ($h \sim k^{-1}$)}
\label{subsec:pollution} 

\subsubsection{Experiments with ImpRAS1 and ImpHRAS} 
 
In this section we consider the discretisation of \eqref{eq:PDE}  with $\eps = 0$ and  $h = \pi/5k$ (i.e. 10 grid points per wavelength). In this case the domain decomposition is load-balanced at about $H = k^{-0.6}$ and so 
we investigated the performance of preconditioned GMRES only for    
$H = k^{-\alpha}$, with  
$\alpha$ in the range $[0.4,0.8]$. We found, for all choices of $\alpha$,
the method HRAS not to be effective (with or without coarse grid solve), and so we focused attention on ImpHRAS and its one-level variant ImpRAS1.

{Sample results for ImpRAS1 (top) and ImpHRAS (bottom) are given in Table \ref{tab:4}.  % ImpHRAS preconditioner   
% % performed best  in the range $[0.4,0.5]$, with sample results
% are given in Table \ref{tab:4} (right hand side).
Here 
$T$ denotes the timing  
for the total solve process, while   $T_{\mathrm{it}}$ denotes the 
time per iteration.   
% For ImpHRAS with $\Hprec = k^{-0.5}$, using the last 6 entries of each column, we see that \# GMRES  grows with about $\mathcal{O}(n^{0.25})$ while  $T$  grows with $\mathcal{O}(n^{1.36})$. 
% with about $O(k^{  })$. % This row also contains the analogous results
% % to $T_{\mathrm{it}}$.
Here the cost of the coarse grid
solve is relatively small and the time per iteration for ImpHRAS is almost
the same as that for ImpRAS1.  
% and  
% and include also time per iteration $T_{\mathrm{it}}$. 
Overall ImpRAS1 is slightly quicker than  ImpHRAS: 
Using the last 6 entries of each column for ImpHRAS with $H = k^{-0.4}$,    
\#GMRES  is  growing with   order $\mathcal{O}(n^{0.18})$, while the total time is growing with order $\mathcal{O}(n^{1.5})$.  
%\ednote{Ivan will check all convergence raates once more} }

% \ednote{add some comments on the tables. Eero can you go higher with
%   $k$? If the iteration numbers are flat we can probably get better
%   rate of growth figures. It does seem that the cost of the direct
%   solver is not proportional to $n$, more line $n^\mu$ with $\mu$
%   given below. }

\begin{table} 
\begin{center} 

{
\begin{tabular}{|c|c||c|c|c|c|c|c|}
\hline 
 &  & \multicolumn{6}{c|}{ImpRAS1}\tabularnewline
\hline 
 &  & \multicolumn{3}{c|}{$H=k^{-0.5}$} & \multicolumn{3}{c|}{$H=k^{-0.4}$}\tabularnewline
\hline 
$k$ & $n$ & $\#GMRES$ & $T$ & $T_{it}$ & $\#GMRES$ & $T$ & $T_{it}$\tabularnewline
\hline 
60 &   9409 &  35 & 6.83 &  0.15 &  20 & 4.67 &  0.16\tabularnewline
80 &  16129 &  39 & 13.01 & 0.27 &  23 & 9.21 &  0.30\tabularnewline
100 &  25921 &  43 & 24.21 & 0.47 &  25 & 18.8 & 0.59\tabularnewline
120 &  35344 &  45 & 37.10 & 0.69 &  29 & 29.50 & 0.83\tabularnewline
140 &  52441 &  49 & 63.85 &  1.12 &  28 & 43.31 & 1.27\tabularnewline
160 &  68121 &  51 & 84.65 & 1.43 &  33 & 67.15 &  1.73\tabularnewline
180 &  82369 &  54 & 113.86 &  1.85 &  32 & 91.01 & 2.43\tabularnewline
200 & 104329 &  57 &  159.67 & 2.47 &  30 & 114.27 &  3.26\tabularnewline
220 & 119716 &  59 & 190.50 &  2.86 &  34 & 160.46 &  4.11\tabularnewline
240 & 141376 &  61 & 249.48 & 3.64 &  35 & 203.30 & 5.12\tabularnewline
260 & 173889 &  66 & 323.79 &  4.43 &  35 & 262.77 & 6.67\tabularnewline
280 & 196249 &  70 & 390.81 &  5.07 &  39 & 354.60 & 8.17\tabularnewline
300 & 227529 &  68 & 459.72 &  6.13 &  38 & 420.12 & 9.98\tabularnewline
\hline 
\end{tabular} }
\quad

 \vspace{0.2cm} 

{ 
\begin{tabular}{|c|c|c|c|c|c|c|c|}
\hline 
\multicolumn{8}{|c|}{ImpHRAS}\tabularnewline
\hline 
\multirow{1}{*}{} &  & \multicolumn{3}{c|}{$H=k^{-0.5}$} & \multicolumn{3}{c|}{$H=k^{-0.4}$}\tabularnewline
\hline 
$k$ & $n$ & $\#GMRES$ & $T$ & $T_{it}$ & $\#GMRES$ & $T$ & $T_{it}$\tabularnewline
\hline 
60 &   9409 &  33 & 5.09 & 0.11 &  21 & 4.36 & 0.14\tabularnewline
80 &  16129 &  40 & 10.87 & 0.22 &  25 & 9.18 & 0.29\tabularnewline
100 &  25921 &  43 & 20.80 & 0.40 &  24 & 17.11 & 0.57\tabularnewline
120 &  35344 &  47 & 34.08 &  0.61 &  29 & 27.99 & 0.79\tabularnewline
140 &  52441 &  52 & 61.25 &  1.01 &  27 & 40.45 &  1.24\tabularnewline
160 &  68121 &  55 & 82.11 &  1.28 &  32 & 63.16 & 1.67\tabularnewline
180 &  82369 &  53 & 103.99 & 1.69 &  32 & 88.40 &  2.37\tabularnewline
200 & 104329 &  56 & 147.72 & 2.29 &  31 & 115.10 &  3.20\tabularnewline
220 & 119716 &  59 & 180.19 &  2.66 &  35 & 161.90 &  4.05\tabularnewline
240 & 141376 &  60 & 233.24 & 3.42 &  35 & 198.54 & 4.99\tabularnewline
260 & 173889 &  64 &  295.08 & 4.05 &  34 & 252.30 &  6.56\tabularnewline
280 & 196249 &  69 & 361.86 & 4.63 &  37 & 332.66 & 8.01\tabularnewline
300 & 227529 &  67 & 430.55 & 5.69 &  37 & 403.04 & 9.76\tabularnewline
\hline 
\end{tabular} }

\end{center} 
\caption{\label{tab:4} Performance of ImpRAS1  (top) and ImpHRAS (bottom) 
with 
$\epsprob = 0$, $\epsprec = k$ and $h = \pi/5k$,   for  $\Hprec = k^{-0.5}, \, k^{-0.4} .$} 
\end{table} 

{In Table \ref{tab:4p} we give preliminary timing results for a parallel implementation
of the ImpRAS1 method. The implementation is in python and is based on numpy and scipy
with the  mpi4py library used for message passing. The problem is run on $P=M^2$ processes, where $M^2$ is
the number of subdomains in the preconditioner. Processes are mapped onto $M$ cluster nodes with $M$
processes running on each node. The column labelled $P$ is the number of processors, which coincides with the number of subdomains. 
The column labelled $n_{loc}$ gives the dimension of the local problem being solved on each processor. 
Note that $n_{loc}$ grows with about $k^{1.2}$ while $P$ grows with about $k^{0.8}$ in this implementation. 
$T$ is the serial time,  $T_{par}$ is the parallel time  and $S = T/T_{par}$.   
Based on the last $6$ entries of the column $T_{par}$, the parallel solve time is growing with about $\mathcal{O}(k^{2.1}) = \mathcal{O}(n^{1.05})$ where $n$ is the system dimension.}

% \begin{table}
%  \begin{center}
%  {\color{red} 
% \begin{tabular}{|c|c|c||c|c|c|c|c|}
% \hline
% $k$ & $P=M^{2}$ & $n_{loc}$ & $\#GMRES$ & $T$ & $T_{par}$ & $S$ & $E$\tabularnewline
% \hline
% \hline
% 60 & 25 & 1444 & 20 & 4.67 & 0.38 & 12.25 & 0.49\tabularnewline
% \hline
% 80 & 36 & 1764 & 23 & 9.21 & 0.51 & 17.97 & 0.50\tabularnewline
% \hline
% 100 & 36 & 2916 & 25 & 18.8 & 1.02 & 18.54 & 0.52\tabularnewline
% \hline
% 120 & 49 & 2916 & 29 & 29.50 & 1.15 & 25.62 & 0.52\tabularnewline
% \hline
% 140 & 49 & 3969 & 28 & 43.31 & 1.62 & 26.66 & 0.54\tabularnewline
% \hline
% 160 & 64 & 3969 & 33 & 67.15 & 1.93 & 34.76 & 0.54\tabularnewline
% \hline
% 180 & 64 & 5041 & 32 & 91.01 & 2.37 & 38.43 & 0.60\tabularnewline
% \hline
% 200 & 64 & 6241 & 30 & 114.27 & 3.05 & 37.43 & 0.58\tabularnewline
% \hline
% 220 & 81 & 6084 & 34 & 160.46 & 3.24 & 49.53 & 0.61\tabularnewline
% \hline
% 240 & 81 & 6889 & 35 & 203.30 & 4.14 & 49.11 & 0.60\tabularnewline
% \hline
% 260 & 81 & 8281 & 35 & 262.77 & 5.34 & 49.23 & 0.61\tabularnewline
% \hline
% 280 & 100 & 8100 & 39 & 354.60 & 5.71 & 62.15 & 0.62\tabularnewline
% \hline
% 300 & 100 & 9025 & 38 & 420.12 & 6.73 & 62.43 & 0.62\tabularnewline
% \hline
% \end{tabular}}
% \end{center} 
% {\color{red} 
% \caption{\label{tab:4p}Parallel performance of ImpRAS1 with $\eps_{prob}=0$,
% $\epsprec=k$ and $h=\pi/5k$, for $\Hprec=k^{-0.4}$. Relative speedup $S$ is  
% shown for comparison of total time $T_{par}$ on
% $P$ processes with serial implementation time $T$.}
% }
% \end{table}

\begin{table}
 \begin{center}
 {
\begin{tabular}{|c|c|c||c|c|c|c|}
\hline
$k$ & $P=M^{2}$ & $n_{loc}$ & $\#GMRES$ & $T$ & $T_{par}$ & $S$ \tabularnewline
\hline
\hline
60 & 25 & 1444 & 20 & 4.67 & 0.38 & 12.25 \tabularnewline
\hline
80 & 36 & 1764 & 23 & 9.21 & 0.51 & 17.97 \tabularnewline
\hline
100 & 36 & 2916 & 25 & 18.8 & 1.02 & 18.54 \tabularnewline
\hline
120 & 49 & 2916 & 29 & 29.50 & 1.15 & 25.62 \tabularnewline
\hline
140 & 49 & 3969 & 28 & 43.31 & 1.62 & 26.66 \tabularnewline
\hline
160 & 64 & 3969 & 33 & 67.15 & 1.93 & 34.76 \tabularnewline
\hline
180 & 64 & 5041 & 32 & 91.01 & 2.37 & 38.43 \tabularnewline
\hline
200 & 64 & 6241 & 30 & 114.27 & 3.05 & 37.43 \tabularnewline
\hline
220 & 81 & 6084 & 34 & 160.46 & 3.24 & 49.53 \tabularnewline
\hline
240 & 81 & 6889 & 35 & 203.30 & 4.14 & 49.11 \tabularnewline
\hline
260 & 81 & 8281 & 35 & 262.77 & 5.34 & 49.23 \tabularnewline
\hline
280 & 100 & 8100 & 39 & 354.60 & 5.71 & 62.15 \tabularnewline
\hline
300 & 100 & 9025 & 38 & 420.12 & 6.73 & 62.43 \tabularnewline
\hline
\end{tabular}}
\end{center} 
{ 
\caption{\label{tab:4p}Parallel performance of ImpRAS1 with $\eps_{prob}=0$,
$\epsprec=k$ and $h=\pi/5k$, for $\Hprec=k^{-0.4}$. Relative speedup $S$ is  
shown for comparison of total time $T_{par}$ on
$P$ processes with serial implementation time $T$.}
}
\end{table}

% \begin{center}
% {\color{blue} 
% \begin{tabular}{|c|c||c|c|c|c|c|}
% \hline 
% $k$ & $P=M^{2}$ & $\#GMRES$ & $T$ & $T_{par}$ & $S_P$ & $E_P$\tabularnewline
% \hline 
% \hline 
% 60  & 25 &  20 & 4.67 & 0.38 & 12.25 & 0.49\tabularnewline
% \hline 
% 80  & 36 &  23 & 9.21 & 0.51 & 17.97 & 0.50\tabularnewline
% \hline 
% 100  & 36 &  25 & 18.8 & 1.02 & 18.54 & 0.52\tabularnewline
% \hline 
% 120  & 49 &  29 & 29.50 & 1.15 & 25.62 & 0.52\tabularnewline
% \hline 
% 140  & 49 &  28 & 43.31 & 1.62 & 26.66 & 0.54\tabularnewline
% \hline 
% 160  & 64 &  33 & 67.15 & 1.93 & 34.76 & 0.54\tabularnewline
% \hline 
% 180  & 64 &  32 & 91.01 & 2.37 & 38.43 & 0.60\tabularnewline
% \hline 
% 200  & 64 &  30 & 114.27 & 3.05 & 37.43 & 0.58\tabularnewline
% \hline 
% 220  & 81 &  34 & 160.46 & 3.24 & 49.53 & 0.61\tabularnewline
% \hline 
% 240  & 81 &  35 & 203.30 & 4.14 & 49.11 & 0.60\tabularnewline
% \hline 
% 260  & 81 &  35 & 262.77 & 5.34 & 49.23 & 0.61\tabularnewline
% \hline 
% 280  & 100 &  39 & 354.60 & 5.71 & 62.15 & 0.62\tabularnewline
% \hline 
% 300  & 100 &  38 & 420.12 & 6.73 & 62.43 & 0.62\tabularnewline
% \hline 
% \end{tabular}
% }

\subsubsection{A multilevel version of ImpRAS1} 

{From Table \ref{tab:4} we see that the  case $H = k^{-0.4}$ provides a  solver 
with remarkably stable  iteration counts, having   
almost no growth with respect to $k$. However (although the coarse 
grid component of the preconditioner can be neglected), the 
local systems to be solved at  each iteration are relatively large, being 
of dimension  $\mathcal{O}((k^{0.6})^2) = \cO(k^{1.2})$. 
We therefore consider inner-outer iterative methods 
where these large  local problems are resolved by an inner 
GMRES preconditioned with an ImpRAS1 preconditioner based on 
decompositon of the local domains of diameter  $k^{-0.4}$ into much smaller 
domains of 
diameter $(k^{-0.4})^2 = k^{-0.8}$. (Such inner-outer methods are also investigated in different ways in \cite{ZeDe:15} and \cite{ShChGr:16}.)  The local problems to be solved then are of dimension $\mathcal{O}((k^{0.2})^2) = \mathcal{O}(k^{0.4})$ and there are 
$\mathcal{O}(k^{1.6})$ of them to solve at each iteration. 

The inclusion  of this method in the present paper is rather tentative, because 
(for the range of $k$ considered), breaking up the local problems of size 
$\mathcal{O}(k^{1.2})$   
into smaller subproblems is not competitive time-wise with the direct solver in  2D.  The times of this multilevel variant are far inferior   to those reported in Table \ref{tab:4}.     
However even though the inner tolerance is set quite large at $0.5$, 
the (outer)  
iteration numbers are remarkably unaffected  
(sample results are given in Table \ref{tab:4a}). 
In this table the outer tolerance is (as before) relative 
residual reduction  of $10^{-6}$. 
Similar results (although slightly inferior) are obtained with $\epsprec = k$, in which case the inner iterations are also  almost identical with those reported 
in Table \ref{tab:4} for ImpRAS1 in the case $\Hprec = k^{-0.4}$. 
   
Since the action of this  preconditioner involves the 
solution of $\mathcal{O}(k^{1.6})$ independent 
local systems of dimension only $\mathcal{O} (k^{0.4})$,  
this method has strong parallel potential and  is also worth investigating  
in 3D, where the direct solvers are less competitive.  

\begin{table}
\begin{center} 
\begin{tabular}{|c|c|c|}
\hline 
$k\backslash\beta$  & 1.2  & 1.6 \tabularnewline
\hline 
100 &  26(6) &  31(4)\tabularnewline
120 &  31(6) &  36(4)\tabularnewline
140 &  29(6) &  35(4)\tabularnewline
160 &  33(7) &  39(5)\tabularnewline
180 &  33(7) &  38(5)\tabularnewline
200 &  32(7) &  39(5)\tabularnewline
220 &  35(8) &  42(5)\tabularnewline
240 &  35(8) &  42(5)\tabularnewline
260 &  34(8) &  42(5)\tabularnewline
280 &  39(9) &  45(6)\tabularnewline
300 &  39(9) &  45(6)\tabularnewline
\hline 
\end{tabular}
\end{center} 
\caption{\label{tab:4a} 
Sample iteration counts for the inner-outer ImpRAS1 preconditioner   
$\epsprob = 0$, $\hprob = \pi/5k$,  $\epsprec = k^{\beta}$,   
$\Hprec = k^{-0.4}$ (for the outer iteration) and $\Hprec = k^{-0.8} $ (for the inner iteration). } 
\end{table}

}

\subsection{Variable wave speed ($h \sim \omega^{-3/2}$)} 
\label{subsec:Variable}
In this subsection we  give some initial results on the performance of our  algorithms when applied to 
problems with variable wave speed. A more  detailed investigation of this problem is one of our 
next priorities and  the discussion here should be regarded as somewhat preliminary. 

Domain 
decomposition methods have  the advantage that the subdomains (and possibly 
the coarse mesh) can be chosen to resolve jumps in the wave speed, if the wave speed 
is geometrically simple enough. At present the variable speed case is not covered by any theory, so this section 
is necessarily experimental. 

We consider  the analogue of the problem \eqref{eq:PDE} with    $k = \omega/c$  
where $\omega$ is the  angular frequency and $c=c(x)$ is the spatially dependent wave speed. 
For the preconditioners we  consider approximate inverses of  problems with variable absorption of the form:    
\begin{equation}\label{eq:PDE1} 
-\Delta u  - (1+ \ri \rho)\left(\frac{\omega}{c}\right)^2 u = f \ , \quad \text{on} \quad \Omega,  \end{equation} 
on a bounded domain $\Omega$ with impedance boundary condition
\begin{equation}\label{eq:ImpBC1} 
\frac{\partial u }{\partial n} - i \left(\frac{\omega}{c}\right) u = g \quad \text{on} \quad \Gamma 
\end{equation}
where 
 $\rho   = \rhoprec \geq  0$ is a  parameter to be chosen.    
Thus when $c$ is constant, and $k := \omega/c$,  the perturbed wavenumber is $ k^2 + \ri \rho k^2$
and so the choice $\eps = k^\beta$ in \eqref{eq:PDE}  corresponds to the choice $\rho = k^{\beta - 2}$ in \eqref{eq:PDE1}.   On the other hand when $c$ is variable, the amount of absorption added is 
proportional to $(\omega/c)^2$ so more absorption is effectively added  
where $c$ is relatively small and less is added  when $c$ is relatively large.  
We do not insert any absorption into the boundary condition \eqref{eq:ImpBC1}. 

We consider a test problem where $\Omega$ is the unit square. An internal square $\Omega_1$ of 
side length $1/3$ is placed  inside  $\Omega$ and the wave speed is taken 
to have value $c^*$ in the 
inner square and  value $1$ in $\Omega_2 : = \Omega \backslash \Omega_1$.  
The square $\Omega_1$  is either placed in the centre of $\Omega$ (this is the case ``discontinuity resolved'', where the coarse grid described below will resolve the interface)  
 or at a position a few fine grid elements  to the north and west of centre, with the distance moved in the directions north and west equal to the size of the overlap of the subdomains.  
In the latter  case the coarse grid passes through the interface (and this is 
called ``discontinuity unresolved'' below).   
We perform experiments with $c^*$ both 
bigger than $1$ and less than $1$ with the latter case expected to be hardest.    

The problem is discretised by a uniform fine grid with $\hprob \sim \omega^{-3/2}$ 
and with the fine grid resolving the interface $\Gamma_{1,2}$ between $\Omega_1$ and $\Omega_2$. 
No absorption is added to the problem to be solved, i.e. $\rhoprob =0$. 

We apply  the inner-outer algorithm as described in \S \ref{subsec:pollution-free} (see Table \ref{tab:3}) for this problem. 
 The outer solver is HRAS with $\Hprec \sim k^{-1}$ while the inner solver is ImpRAS1 with $\Hprec \sim k^{-1/2}$. For both inner and outer solvers we set   
$\rho_{prec} = \omega^{\beta-2}$. 
In all cases generous overlap is used and the RAS domains are determined by the coarse grid as described in the introductory paragraphs to this section. 
  
The coarse grid for the outer solve consists of uniform triangles of diameter $ \sim k^{-1}$ which are chosen to  resolve the square $\Omega_1$ when it is placed in the centre, and do not resolve it when the square is moved. Numerical results,  comparing  the cases $c^* =  1.5, 1,  0.66$ 
are given in Table \ref{tab:6a} \ref{tab:6b}, \ref{tab:6c}. In each row, for each value of $\beta$, 
 the three figures indicate the number of outer HRAS iterations, the number of inner ImpRAS1 iterations (in brackets) and the total time on a serial machine.  
{The outer tolerance is set at $10^{-6}$ while the inner tolerance
is set at { $0.5$}. 

The times for $\beta = 1.6$  grow with about $\mathcal{O}(n^{1.4})$ in the case $c^* = 1.5$ and $c^* = 1$ (rather similar to the performance observed in Table \ref{tab:3}). The actual times in the case $c^* = 0.66$ are considerably worse (which is to be expected as smaller $c*$ implies larger effective frequency on that domain. But the rate of growth of time  with $n$ is not affected very much, being about $\mathcal{O}(n^{1.5})$ in 
Table \ref{tab:6c}. The case $c^* = 1.5 $ seems a little easier to solve than the case $c^*= 1$. There is not much difference in any case between the resolved and the unresolved cases.     

\begin{center}
\begin{table}
\begin{center}
\begin{tabular}{|c|cc|cc|cc|cc|}
\hline 
\multicolumn{9}{|c|}{$c^{*}=1.5$, discontinuity resolved}\tabularnewline
\hline 
$\omega\backslash\beta$  & \textbf{1.0 } &  & \textbf{1.2 } &  & \textbf{1.6 } &  & \textbf{1.8} & \tabularnewline
\hline 
\multirow{1}{*}{\textbf{10} } &  19(1) &    0.71 &  19(1) &    0.55 &  20(1) &    0.53 &  21(1) &    0.54\tabularnewline
\multirow{1}{*}{\textbf{20} } &  20(2) &    3.25 &  20(2) &    3.22 &  27(1) &    3.65 &  30(1) &    3.84\tabularnewline
\multirow{1}{*}{\textbf{40} } &  22(3) &   50.09 &  23(3) &   50.55 &  29(2) &   54.04 &  44(1) &   62.99\tabularnewline
\multirow{1}{*}{\textbf{60} } &  25(4) &  356.71 &  26(4) &  358.10 &  35(2) &  381.06 &  57(1) &  445.19\tabularnewline
\multirow{1}{*}{\textbf{80 }} &  29(5) & 1244.13 &  29(4) & 1240.80 &  40(2) & 1394.72 &  66(1) & 1606.64\tabularnewline
\multirow{1}{*}{\textbf{100 }} &  35(6) & 3479.95 &  35(5) & 3697.02 &  45(2) & 3820.97 &  78(1) & 4309.29\tabularnewline
\hline 
\multicolumn{9}{|c|}{$c^{*}=1.5$, discontinuity unresolved}\tabularnewline
\hline 
$\omega\backslash\beta$  & \textbf{1.0 } &  & \textbf{1.2 } &  & \textbf{1.6 } &  & \textbf{1.8} & \tabularnewline
\hline 
\multirow{1}{*}{\textbf{10} } &  18(1) &    0.70 &  19(1) &    0.56 &  20(1) &    0.53 &  21(1) &    0.54\tabularnewline
\multirow{1}{*}{\textbf{20 }} &  20(2) &    3.26 &  20(2) &    3.25 &  27(1) &    3.65 &  30(1) &    3.87\tabularnewline
\multirow{1}{*}{\textbf{40 }} &  22(3) &   50.80 &  23(3) &   51.30 &  29(2) &   54.56 &  44(1) &   63.76\tabularnewline
\multirow{1}{*}{\textbf{60 }} &  25(4) &  363.19 &  26(4) &  364.96 &  35(2) &  387.40 &  58(1) &  454.04\tabularnewline
\multirow{1}{*}{\textbf{80 }} &  30(5) & 1273.11 &  30(4) & 1347.66 &  40(2) & 1417.61 &  66(1) & 1623.74\tabularnewline
\multirow{1}{*}{\textbf{100 }} &  35(6) & 3545.44 &  35(5) & 3541.62 &  45(2) & 3660.37 &  78(1) & 4042.62\tabularnewline
\hline 
\end{tabular}
\caption{\label{tab:6a} 
Performance of the inner-outer algorithm described in \S \ref{subsec:Variable}.
Discontinuous wave speed, $c* = 1.5$. 
  } 
\par\end{center}
\end{table} 
\end{center} 

\begin{table} 
\begin{center}
\begin{tabular}{|c|cc|cc|cc|cc|}
\hline 
\multicolumn{9}{|c|}{$c^{*}=1.0$}\tabularnewline
\hline 
$\omega\backslash\beta$  & \textbf{1.0 } &  & \textbf{1.2 } &  & \textbf{1.6 } &  & \textbf{1.8} & \tabularnewline
\hline 
\multirow{1}{*}{\textbf{10} } &  18(1) &    0.70 &  18(1) &    0.54 &  19(1) &    0.51 &  21(1) &    0.54\tabularnewline
\multirow{1}{*}{\textbf{20} } &  19(2) &    3.12 &  19(2) &    3.16 &  25(1) &    3.43 &  29(1) &    3.73\tabularnewline
\multirow{1}{*}{\textbf{40} } &  22(3) &   48.76 &  22(3) &   48.58 &  28(2) &   51.79 &  45(1) &   62.22\tabularnewline
\multirow{1}{*}{\textbf{60} } &  28(5) &  353.26 &  28(4) &  352.74 &  35(2) &  368.99 &  56(1) &  429.53\tabularnewline
\multirow{1}{*}{\textbf{80 }} &  36(5) & 1253.44 &  35(5) & 1244.01 &  42(2) & 1361.78 &  66(1) & 1476.53\tabularnewline
\multirow{1}{*}{\textbf{100 }} &  45(7) & 3487.02 &  44(6) & 3693.13 &  49(2) & 3728.06 &  79(1) & 4179.60\tabularnewline
\hline 
\end{tabular}
\par\end{center}
\caption{\label{tab:6b} Performance of the inner-outer algorithm described in \S \ref{subsec:Variable}.
 Continuous wave speed $c* = 1$}
\end{table} 

\begin{table}
\begin{center}
\begin{tabular}{|c|cc|cc|cc|cc|}
\hline 
\multicolumn{9}{|c|}{$c^{*}=0.66$, discontinuity resolved}\tabularnewline
\hline 
$\omega\backslash\beta$  & \textbf{1.0 } &  & \textbf{1.2 } &  & \textbf{1.6 } &  & \textbf{1.8} & \tabularnewline
\hline 
\multirow{1}{*}{\textbf{10} } &  19(1) &    0.73 &  20(1) &    0.58 &  21(1) &    0.54 &  23(1) &    0.58\tabularnewline
\multirow{1}{*}{\textbf{20} } &  22(2) &    3.38 &  22(2) &    3.43 &  28(1) &    3.68 &  33(1) &    4.03\tabularnewline
\multirow{1}{*}{\textbf{40} } &  31(4) &   55.03 &  32(3) &   55.22 &  37(2) &   57.46 &  54(1) &   68.06\tabularnewline
\multirow{1}{*}{\textbf{60} } &  48(5) &  418.78 &  48(4) &  415.63 &  54(2) &  426.52 &  79(1) &  502.58\tabularnewline
\multirow{1}{*}{\textbf{80 }} &  85(7) & 1709.73 &  78(5) & 1628.70 &  74(2) & 1630.55 & 108(1) & 1925.38\tabularnewline
\multirow{1}{*}{\textbf{100 }} & 124(8) & 4881.62 & 115(7) & 4853.22 &  93(2) & 4448.73 & 134(1) & 5151.77\tabularnewline
\hline 
\multicolumn{9}{|c|}{$c^{*}=0.66$, discontinuity unresolved}\tabularnewline
\hline 
$\omega\backslash\beta$  & \textbf{1.0 } &  & \textbf{1.2 } &  & \textbf{1.6 } &  & \textbf{1.8} & \tabularnewline
\hline 
\multirow{1}{*}{\textbf{10} } &  19(1) &    0.72 &  19(1) &    0.60 &  21(1) &    0.54 &  23(1) &    0.58\tabularnewline
\multirow{1}{*}{\textbf{20 }} &  23(2) &    3.54 &  23(2) &    3.55 &  29(1) &    3.84 &  34(1) &    4.19\tabularnewline
\multirow{1}{*}{\textbf{40 }} &  32(4) &   58.89 &  32(3) &   58.45 &  38(2) &   61.33 &  55(1) &   71.74\tabularnewline
\multirow{1}{*}{\textbf{60 }} &  49(5) &  450.16 &  49(4) &  448.57 &  55(2) &  458.60 &  80(1) &  533.56\tabularnewline
\multirow{1}{*}{\textbf{80 }} &  85(7) & 1820.58 &  79(5) & 1826.07 &  77(2) & 1767.87 & 109(1) & 2041.53\tabularnewline
\multirow{1}{*}{\textbf{100 }} & 123(8) & 5076.60 & 116(6) & 5016.11 &  96(2) & 4567.90 & 135(1) & 5323.54\tabularnewline
\hline 
\end{tabular}
\par\end{center}
\caption{\label{tab:6c}  Performance of the inner-outer algorithm described in \S \ref{subsec:Variable}.
 Discontinuous wave speed, $c* = 0.66$} 
\end{table}

\

%\subsubsection*{ImpHRAS - the case c{*}=1.0:}

{
\section{Summary} 

%\ednote{Euan: could you look this over please?} 
%\ednote{Ivan to write what is the RHS in the experiments} 
In this paper we considered the  construction of preconditioners for 
the Helmholtz equation (without or with  absorption)  by  using  
domain decomposition methods applied  
to the corresponding  problem 
with absorption. 

These methods are  related to the shifted Laplacian multigrid methods, but the relative simplicity of the method considered here   permits  
rigorous analysis of the convergence of GMRES through estimates of the field of values of the preconditioned problem.  The  flexibility of the domain decomposition approach also allows for 
the insertion of sub-solvers which are appropriate for high frequency Helmholtz problems, such as replacing Dirichlet local problems with impedance (or PML) local problems.      

For the analysis, two theoretical subproblems  are  identified:  (i) What range of  absorption {\em is permitted},  so that   the problem 
with absorption remains an optimal  preconditioner for the problem 
without absorption? 
 and (ii) What range of absorption {\em is needed} so that  the 
domain decomposition method performs  optimally as 
a preconditioner for the problem with absorption?    

The ranges  
that  result from studying  problems (i) and (ii) separately   
have been analysed,  and this analysis  is reviewed in the paper (\S \S \ref{sec:Intro}, \ref{sec:Main}). Since these ranges are  disjoint,  the  best methods are obtained by using 
a combination of insight provided by the rigorous analysis and  
by numerical experimentation. The best methods involve careful tuning of 
the absorption parameter, the choice of coarse grid and the choice of boundary condition on the subdomains.    

Of those methods studied, the best 
(in terms of computation time on a serial machine) differ,  depending on 
the level of resolution of the underlying finite element grid. For  
problems with  constant wave speed and with mesh diameter $h \sim k^{-3/2}$ (so chosen to resolve the pollution effect), a multilevel method with serial time complexity 
 $\mathcal{O}(n^\alpha)$ with $\alpha \in [1.3,1.4]$ is presented, where $n \sim k^3$ is the dimension of the system being solved (\S\ref{subsec:pollution-free}). In this method a two level  
preconditioner with a  fairly fine coarse grid is used, and the coarse grid problem is resolved by an inner iteration with a further one-level preconditioner with impedance local solves.  

For discretisations involving a fixed number of grid points per wavelength, similar time complexity is    achieved by highly parallelisable 
one-level methods using impedance local solves on relatively large subdomains. 

We also illustrate the method when it is applied to a model problem 
with jumping wave speed (\S \ref{subsec:Variable}). A preliminary parallel experiment is  also given. 
}

\end{document}